\documentclass[12pt,a4paper]{article}
\usepackage[utf8]{inputenc}
\usepackage{graphicx} % required for inserting image 
\usepackage{array,mathtools}
\usepackage[square,sort,comma,numbers]{natbib}
\usepackage{amsmath, amssymb}
\usepackage{svg}
\usepackage{url}
\usepackage{subcaption}
\usepackage{booktabs}
\usepackage{hyperref} % hyperlinks 
\hypersetup{
    colorlinks=true,  % Enable colored links
    linkcolor=black,   % Set the color for internal links (e.g., table of contents)
    urlcolor=black,    % Set the color for external links (e.g., URLs)
    citecolor=red
    }

\newenvironment{keywords}{%
\renewcommand{\section}[2]{}% Ignore \section within keywords
\small\textbf{Keywords:}}{%
}

\graphicspath{{./Graphs/}}

\DeclareMathOperator{\sign}{sign}
\DeclareMathOperator{\Si}{Sigm}% 1 in latex 2 in pdf
\DeclareMathOperator{\newS}{S}% 1 in latex 2 in pdf

\DeclareMathOperator{\rhs}{Rhs}
\DeclareMathOperator{\e}{e}

\newcommand{\tran}{\mathsf{T}}
\newcommand{\R}{\mathbb{R}}
\newcolumntype{P}[1]{>{\raggedright\arraybackslash}p{#1}}

\newcommand{\hmdonotuse}[1]{}

\title{Alpha-Delta Transitions in Cortical Rhythms as grazing bifurcations}
\author{Huda Mahdi 
\thanks{Department of Mathematics and Statistics, University of Exeter, Living Systems Institute, Exeter EX44QD, UK.  \protect\url{mailto: hm672@exeter.ac.uk} ,
\protect\url{https://orcid.org/0009-0004-7318-7096}}
\and  Jan Sieber
\thanks{Department of Mathematics and Statistics, University of Exeter, Harrison Building, Exeter EX4 4QF, UK.  \protect\url{mailto: j.sieber@exeter.ac.uk} ,
\protect\url{https://orcid.org/0000-0002-9558-1324}}
\and   Krasimira Tsaneva-Atanasova
\thanks{Department of Mathematics and Statistics and EPSRC Hub for Quantitative Modelling in Healthcare, Living Systems Institute, University of Exeter, Exeter, EX44QJ,
UK.  \protect\url{mailto: k.tsaneva-atanasova@exeter.ac.uk} ,
\protect\url{https://orcid.org/0000-0002-6294-7051}}
}

\date{\today}

\date{\today}

\date{April 2024}

\begin{document}

\maketitle
%JS word limit abstract =200 words => have to shorten slightly
\begin{abstract}

The Jansen-Rit model of a cortical column in the cerebral cortex is widely used to simulate spontaneous brain activity (EEG) and event-related potentials. It couples a pyramidal cell population with two interneuron populations, of which one is fast and excitatory and the other slow and inhibitory. 

Our paper studies the transition between alpha and delta oscillations produced by the model. Delta oscillations are slower than alpha oscillations and have a more complex relaxation-type time profile. In the context of neuronal population activation dynamics, a small threshold means that neurons begin to activate with small input or stimulus, indicating high sensitivity to incoming signals. A steep slope signifies that activation increases sharply as input crosses the threshold.
Accordingly in the model the excitatory activation thresholds are small and the slopes are steep. Hence, a singular limit replacing the excitatory activation function with all-or-nothing switches, eg. a Heaviside function, is appropriate. In this limit we identify the transition between alpha and delta oscillations as a discontinuity-induced grazing bifurcation. At the grazing the minimum of the pyramidal-cell output equals the threshold for switching off the excitatory interneuron population, leading to a collapse in excitatory feedback.
%We identify the boundaries that distinguish these two rhythms employing bifurcation analysis and approximations of neuronal firing rates using a Heaviside step function. 
%Our findings reveal that delta activity is driven by excitatory interneurons, activated when pyramidal cells receive high-amplitude excitatory feedback, while alpha activity occurs under low-amplitude feedback conditions. Specifically, we show that when excitatory input surpasses a critical threshold, excitatory interneurons remain continuously active, while sub-threshold input results in switch-like behaviour, characterised by intermittent activation. The interaction between this switch-like behaviour and the refractory period is key to generating distinct rhythms.

%Physiologically alpha activity is linked to relaxed wakefulness and cognitive engagement, while delta activity is associated with deep sleep and restorative processes. Therefore, understanding these transitions enhances our knowledge of how neural circuits regulate brain states and may provide insights into disruptions that affect sleep and cognitive functions.

\end{abstract}


\section{Introduction}

The Jansen–Rit model \citep{jansen1993neurophysiologically,jansen1995electroencephalogram} is a well-established neural mass model of the activity of a local cortical circuit in the human brain. 
The model builds on the earlier work of da~Silva~et~al.\ \citep{LopesdaSilva1974,LopesdaSilva1976}, \citet{vanRotterdam1982}  and \citet{freeman1987simulation}.  
They developed mathematical models to simulate spontaneous brain activity, which can be recorded non-invasively, to investigate the mechanisms behind specific types of electroencephalogram (EEG) field potentials. Their focus was on understanding the origins of alpha-like activity – a rhythmic EEG pattern with a frequency range of 8-12 Hz, most prominent during restful states with closed eyes. 

Neuronal activity detected by EEG results from the combined excitatory and inhibitory postsynaptic potentials generated by large groups of neurons, such as cortical columns in the cerebral cortex, firing at the same time. \citet{jansen1993neurophysiologically} and  \citet{jansen1995electroencephalogram} extended Lopes da Silva's lumped parameter model \citep{LopesdaSilva1974,LopesdaSilva1976} by incorporating an excitatory feedback loop from  local interneurons in addition to the original populations of inhibitory interneurons and excitatory pyramidal cells in order to investigate the significance of excitatory connections over long distances. Although the Jansen and Rit model \citep{jansen1995electroencephalogram}  as illustrated in Figure \ref{fig:cir} is a simplification of the complexity of neural  connections in cortical regions of the brain, it is able to generate a range of wave forms and rhythms resembling EEG recordings. Accordingly it has been extensively employed to simulate brain rhythmic activity  recorded by EEG  (see \citep{david2003neural,stefanovski2019linking} and  references therein).   
The Jansen-Rit model has been extended further by \citet{wendling2000relevance} to study epileptic-like oscillations by the addition of a fast inhibitory interneuron population. This model has since gained significant attention in studies of epileptic seizures \citep{wendling2000relevance,10.1007/978-3-031-21484-4_25}. Furthermore, the Jansen-Rit model has been used to simulate event-related potentials (ERPs) by applying pulse-like inputs, allowing for the replication and analysis of EEG responses to external stimuli. Notably, the interaction between cortical columns has been found to play a key role in generating visual evoked potentials (VEPs) \citep{jansen1993neurophysiologically}.

The rhythms generated by the Jansen-Rit model can be associated with different behaviours, levels of excitability, and states of consciousness. Using the parameter settings suggested by \citet{jansen1995electroencephalogram}, the model produces oscillations around 10 Hz corresponding to alpha-like activity as described by \citet{grimbert2006bifurcation}. 
Delta rhythms have frequencies between 0.5 and 4 Hz detected during deep stages of non-REM sleep (particularly stages 3 and 4), also known as slow-wave sleep. During these stages, the brain exhibits large-amplitude, low-frequency delta activity. A recent experimental study reports an alpha/delta switch in the prefrontal cortex regulating the shift between positive and negative emotional states \citep{burgdorf2023prefrontal}. Moreover, another very recent study by \citet{Brudzynski2024} identifies a transition between positive (associated with alpha rhythms) and negative (associated with delta rhythms) emotions controlled by the arousal system. At the level of neuronal circuits, the alpha rhythm is linked to synaptic long-term potentiation (LTP) in the cortex, while the delta rhythm is associated with synaptic depotentiation (LTD) in the same region. Therefore, understanding what controls the transitions between alpha and delta oscillations in neural circuits enhances our knowledge of how neural circuits regulate brain states and may provide insights into disruptions that affect sleep and cognitive functions

\citet{grimbert2006bifurcation} analysed the oscillatory behaviour of the Jansen-Rit model using numerical bifurcation analysis. In their work, the external input signal $p$ (see Figure \ref{fig:cir} below) was considered as a bifurcation parameter. The behaviour of the model under variation in $p$ was investigated where all other model parameters were kept fixed and equal to what was originally proposed by Jansen and Rit (see Table \ref{tab:paramter model} below). 
They found that the qualitative behaviour of the system changes from steady state to oscillatory behaviour (appearance or disappearance limit cycle) via a Hopf bifurcation. The authors  also pointed out that the system displays distinct phenomena such as bistability, limit cycles and chaos. Furthermore, \citet{touboul2009codimension} carried out bifurcation analysis for a non-dimensionalised Jansen-Rit model and the extension proposed by \citet{wendling2000relevance} in several system parameters, detecting and tracking bifurcations of up to codimension $3$. 

The transition between alpha and delta oscillations is linked to a notable feature exhibited by the Jansen-Rit neural mass model: a sharp qualitative change in the oscillations' time profiles and frequencies occurs over a small parameter range without change of stability and no (or little) hysteresis. Hence, the underlying mathematical mechanism cannot be a generic bifurcation as found in textbooks \citep{K04,guckenheimer2013nonlinear}. Recent work by \citet{forrester2020role}, which developed an analysis of a large network of interacting neural populations of Jansen-Rit  cortical column  models in the absence of noise identified the alpha-delta transition as a ``\emph{false bifurcation}''. They used an arbitrary geometric feature of the periodic orbit (an inflection point) and  associated it with the qualitative changes of  the periodic orbits. The feature of the orbits could then be expressed as  a zero problem and tracked in the two-parameter space $(A, B)$, where $A$ and $B$ are the gain (strength) of excitatory and inhibitory responses, respectively. This is similar to the approach proposed by \citet{rodrigues2010method}, who tracked qualitative changes in orbit shapes in a
EEG model of absence seizures (such as  \citet{marten2008onset}).
study of alpha and delta frequency oscillations based on different parameter settings has been performed by  \citet{ahmadizadeh2018bifurcation}. They showed that a model of two coupled cortical columns can produce delta activity as the gain strength between the two cortical columns is varied (see Fig.12T1 of \citet{ahmadizadeh2018bifurcation} for delta-like activity). Their model can also produce alpha-like activity, where the frequency of oscillation lies in alpha frequency band but the amplitude changes.   

Our paper starts from the observation that after non-dimensionalisation (as done by  \citet{touboul2009codimension} the activation function for excitatory inputs has steep slopes and small thresholds. The small threshold models that neurons begin to activate with small input or stimulus, indicating high sensitivity to incoming signals. A steep slope signifies that once the activation begins, it escalates rapidly, with the neurons' response intensity increasing sharply as input crosses the threshold. Hence, it makes sense to introduce  a small parameter $\epsilon$ equal to a quarter of the inverse of the activation slope. The singular limit $\epsilon\to0$ replaces the activation functions by all-or-nothing switches represented mathematically by Heaviside functions. In this singular limit we can locate most bifurcations of equilibria with simple formulas. We can also identify the underlying mechanism for the alpha-delta transition. During alpha-type oscillations the minimum of the potential of the pyramidal cells always stays above the threshold for switching off the excitatory feedback (see Figure~\ref{figsamllepstime}). At the transition the minimum then ``touches'' this threshold. If the potential falls below the threshold even briefly the excitatory feedback collapses leading to all potentials dropping to zero. This type of ``touching of a threshold'' is a typical bifurcation in a system with discontinuous switches, called a \emph{grazing bifurcation} \citep{bernardo2008piecewise,kuznetsov2003one}. Our Figure \ref{fig5} shows how this grazing bifurcation is an accurate approximation of the boundary between alpha- and delta-type oscillations for the value $\epsilon\approx 1/40$ corresponding to the parameters in the original Jansen-Rit model.

The paper is organised as follows. We present the Jansen-Rit model in Section \ref{sec:observation}, together with a first numerical bifurcation analysis in one parameter (feedback strength $A$ between excitatory populations). This analysis shows the two types of oscillation (alpha and delta type) and the sharp transition between them. Section~\ref{sec:singpert:JR} first non-dimensionalises the model, identifies a small parameter $\epsilon$ and different ways to take the singular limit $\epsilon\to0$. In Subsection~\ref{subsec:Equilibria in the  limit} we derive explicit expressions for the location of equilibria in the limiting piecewise linear system. In Subsection~\ref{sec:Piece-wise periodic orbits} we derive an algebraic system of equations for the periodic orbits of alpha-type, which are piecewise exponentials and for which one of the components is constant. Using this algebraic system we detect and track the grazing bifurcation, resulting in Figure~\ref{fig5}. As our singular perturbation analysis is only partial, we discuss open questions in Section~\ref{sec:discussion}. 

\section{The Jansen-Rit model for a single cortical column}
\label{sec:observation}
\begin{figure}[ht]
    \centering
 \begin{subfigure}{0.45\textwidth}
        \centering
        \includegraphics[width=0.8\textwidth]{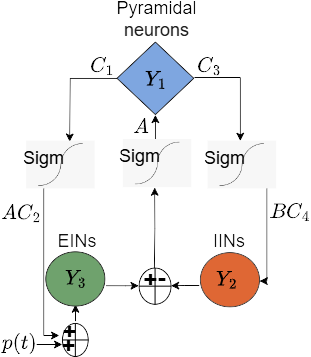}
        \caption{Block diagram Jansen-Rit model \eqref{jansenrit}}
        \label{fig:cir}
    \end{subfigure}
\hfill
    \begin{subfigure}{0.5\textwidth}
        \centering
        \includegraphics[width=\textwidth]{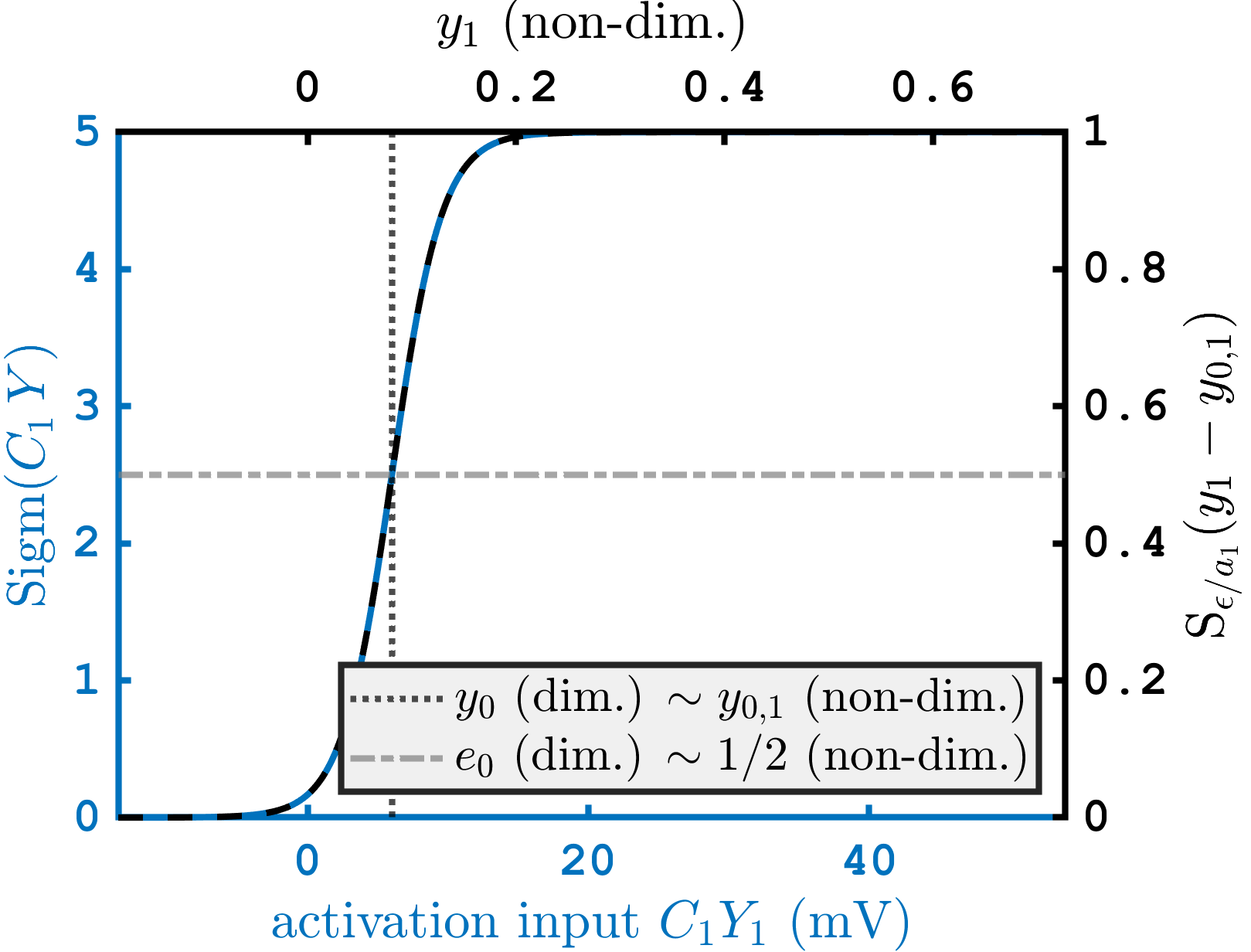}
        %dimnondimsigm
        \caption{Profile of $\newS$ (non-dim.) and $\Si$ (dim.)}
         \label{fig:sigmwithdiffeps}
    \end{subfigure}
    \caption{(a) Interactions between $3$ neuronal populations in local circuit of  a single cortical column in cerebral cortex of brain modelled by \eqref{jansenrit}: pyramidal cells ($Y_1$) and excitatory ($Y_3$, EINs) and inhibitory ($Y_2$, IINs) interneurons. (b) Blue axis: Sigmoid profile of $\Si$, given in \eqref{orginsigm}, with 
    slope $r=0.56$ at threshold $y_0=6$\,mV (vertical dotted line), at which half maximum $e_0$ of $\Si$ is achieved. (b) Black axis:  dimensionless sigmoid $\newS$, given in \eqref{newsigwitheps}, for 
    maximal slope value$1/(4\epsilon)=0.25/0.024\approx10.5$ at activation threshold $y_{0,1}=0.08$ (non-dimensionalised threshold $y_{0,1}$ for excitatory populations indicated by black vertical dotted line).
    See Table~\ref{tab:paramter model} and Table~\ref{tab:nondimparameters} for other parameters of $\Si$ and $\newS$. 
   }
\end{figure}
Figure~\ref{fig:cir} shows the structure of the Jansen-Rit model for a single cortical column. The model assumes that the cortical column contains two populations of interneurons, one excitatory and one inhibitory, and a population of excitatory pyramidal cells. The model contains two blocks for each population, a linear dynamic block, modelled as a second-order differential equation with the population's average \emph{postsynaptic potential} (PSP) $Y_i$ as its output and an average pulse density of action potentials as its input. In Fig.~\ref{fig:cir} these are the blocks with label $Y_i$: $Y_1$ is the PSP of the pyramidal cell population (blue diamond in Fig.~\ref{fig:cir}), $Y_2$ is the PSP of the inhibitory interneuron population (red ellipse in Fig.~\ref{fig:cir}), and $Y_3$ is the PSP of the excitatory interneuron population (green ellipse in Fig.~\ref{fig:cir}). The other block for each population is a nonlinear sigmoid-type activation function (the blocks labelled ``Sigm'' in Fig.~\ref{fig:cir}) from the PSP into an average pulse density of action potentials.

The arrows in Figure  \ref{fig:cir} show a positive feedback loop between the pyramidal neuron population ($Y_1$) and the excitatory interneuron population ($Y_3$), where both connections are excitatory;  and a negative feedback loop between the pyramidal cells and the inhibitory interneurons ($Y_2$), where the connection back from the inhibitory neurons to the pyramidal cells is inhibitory. The model incorporates an external excitatory input, labelled $p(t)$ (an average pulse density) representing signalling from other neuronal populations external to the column (Fig.~\ref{fig:cir}).
The  resulting ODE system corresponding to the schematic diagram Fig.~\ref{fig:cir} is as follows (note that each equation is of order two):
\begin{align} \label{jansenrit} 
\begin{split} 
%\nonumber 
%\dot Y_{1}(t)&=Y_{4}(t)\\
 %\nonumber 
%\dot Y_{2}(t)&=Y_{5}(t)\\
%\nonumber
%\dot Y_3(t)&= Y_6(t)\\
 %\nonumber
 \ddot Y_{1}(t)&=Aa\{\Si(Y_{3}(t)-Y_{2}(t))\}-2a\dot Y_{1}(t)-a^2Y_{1}(t), \\ 
 %\nonumber
\ddot Y_{2}(t)&=Bb\{C_{4}\Si(C_{3}Y_{1}(t)\}-2b\dot Y_{2}(t)-b^2Y_{2},\\
 %\nonumber
 \ddot Y_3(t)&=Aa\{p(t)+C_2\Si(C_{1}Y_1(t))\}-2a\dot Y_3(t)-a^2Y_3(t)\mbox{,\quad where}
 \end{split}\\
%\end{align}
%\begin{align} 
	{\Si(Y)}&=\frac{2e_{0}}{1+\mathrm{e}^{r(y_{0}-Y)}}\label{orginsigm}
 \end{align}
is a nonlinear sigmoid activation function, converting local field postsynaptic potential into firing rate, as shown in Fig.~\ref{fig:sigmwithdiffeps}. Its argument $y$  is the average  postsynaptic potential  %membrane potential state 
in $\text{mV}$ of the neural population. The threshold $y_0$ determines the input level at which the firing rate is half of its maximum $2{e_0}$, and $r$ is the slope of  $\Si$ at $y_0$.

We restrict our study mostly to the case without external input signal by setting the  input $p(t)$ to zero, i.e. $p(t)=0$. In previous studies the values of $p(t)$ has been varied from 120 Hz to 320 Hz as proposed by \citet{jansen1995electroencephalogram}.
For example,  \citet{grimbert2006bifurcation} 
have performed bifurcation analysis of \eqref{jansenrit} %the Jansen and Rit model 
varying 
the input $p(t)$.

\begin{table}[t]
\caption{\label{tab:paramter model}Quantities in Jansen-Rit model, \eqref{jansenrit}, and their default values and dimensions\citep{jansen1995electroencephalogram}.}
    \centering
    \begin{tabular}{lp{0.5\linewidth}P{0.2\linewidth}}
        \toprule % Thicker line at the top
        Parameter & Description & Value/unit \\ 
        \midrule % Line between header and content
        $Y_i(t),i=1,3$&  excitatory postsynaptic potential EPSP & $\text{mV}$ \\
         $Y_2(t)$& inhibitory postsynaptic potential IPSP & $\text{mV}$ \\
        $2e_0$  &asymptotic
maximum of $\Si$  & $5\,\text{s}^{-1}$   \\
       $y_0$ & switching threshold of $\Si$ ($\Si(y_0)=e_0$)  & $6\, [\text{mV}]^{-1}$  \\
       $r$  &  slope of $\Si$ at $y=y_0$  ($\Si'(y_0)$)& $0.56\, [\text{mV}]^{-1}$\\
       $A$ & maximal amplitude of excitatory PSP
       &$3.25\, \text{mV}$\\
       %, 
        $B$ & maximal amplitude of inhibitory PSP  & $22\, \text{mV}$\\
        $a$ & decay rate in excitatory populations& $100\, \text{s}^{-1}$\\
        $b$& decay rate in inhibitory population & $50\, \text{s}^{-1}$\\
$\alpha_i$, $i=1,\ldots,4$& average number of synaptic connections  &$\alpha_1=1,\alpha_2=0.8,\alpha_3=0.25$, $\alpha_4=0.25$\\
$C_i$, $i=1,\ldots,4$&connectivity strength &$C_1=135\alpha_1,\ C_2=\alpha_2C_1,\ C_3=\alpha_3C_1,\ C_4=\alpha_4C_1$\\
$p$& External signal&0\\
        \bottomrule % Thicker line at the bottom
    \end{tabular}
\end{table}

The quantity $Y_3(t)-Y_2(t)$ can be related to experiments, because it is proportional to the signals obtained via EEG recordings corresponding to the average local field potential generated by the underlying neuronal populations in the cortical circuits \citep{jansen1993neurophysiologically}. In the cortex, pyramidal neurons project their apical dendrites into the superficial layers where the post-synaptic potentials are summed, thereby making up the core of the EEG.
The interpretation of all model quantities and their numerical values and dimensions are given in Table \ref{tab:paramter model}. The values are set to those proposed by \citet{jansen1995electroencephalogram}.

\begin{figure}[ht] 
\centering
\includegraphics[width=\textwidth]{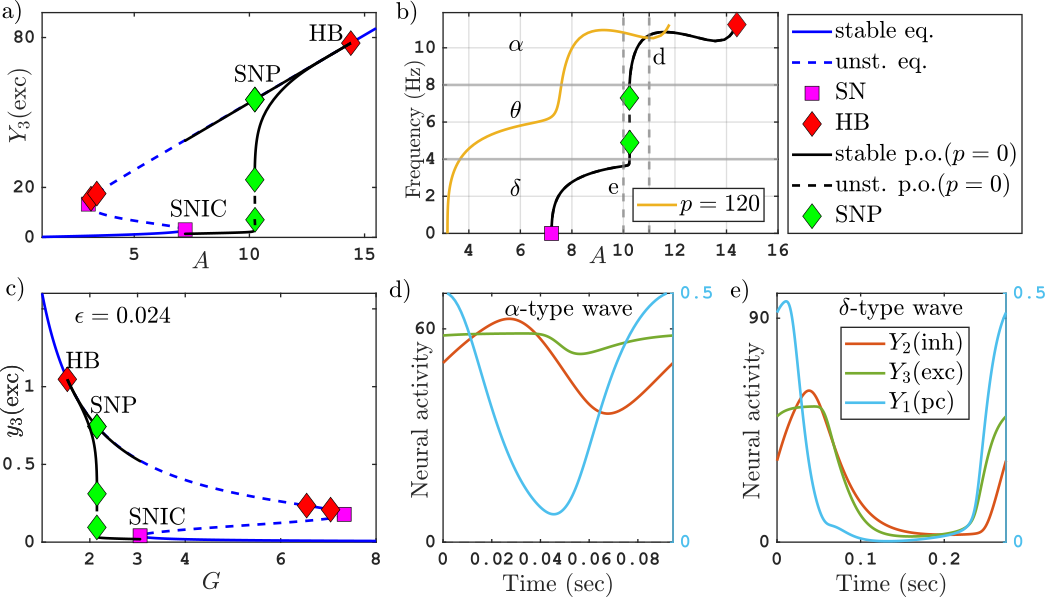} 
\caption{ (a) Bifurcation diagram of Jansen-Rit model \eqref{jansenrit} for varying ${A}$; 
(b)  frequency of oscillations with input ${p=0}$ (black) and input ${p=120}$ (yellow) for varying ${A}$, indicating alpha, theta and delta frequency ranges. (c) Same bifurcation diagram as (a) but using nondimensionalised quantities $G=B/A$ and $y_3$. (d),(e) Time profiles of oscillations in alpha (${A=11}$) and delta (${A=10}$) rhythm regimes corresponding to vertical dashed lines in panel b. See Tables \ref{tab:paramter model} and \ref{tab:nondimparameters} for parameters. Computation performed with \textsc{coco} \citep{dankowicz2013recipes}.} 
\label{dimandnondim}
\end{figure}
\paragraph{Sharp transitions from alpha  to delta activity}
Figure~\ref{dimandnondim}a shows the one-parameter bifurcation diagram of the  Jansen and Rit  model \eqref{jansenrit} for varying excitability scaling factor ${A}$, which determines maximal amplitude  of excitatory postsynaptic potentials (EPSP)  of excitatory populations (pyramidal cells $Y_1$) and exictatory interneurons ($Y_3$).
%which controls the connection strength between? the exictatory populations (pyramidal cells $Y_1$) and exictatory interneurons ($Y_3$). 
All other parameters are fixed as listed in Table~\ref{tab:paramter model}. The $y$-axis shows the PSP of the excitatory interneurons ($Y_3$). 
%
%\hmq{correct the meaning  feedback strenghth of exictatory  }
%
A branch of equilibria (blue) folds back and forth in
two saddle node bifurcations (purple squares). Solid curves are stable, dashed curves are unstable parts of the branch. Along the upper part of the branch the equilibria change stability in three Hopf bifurcation (red diamonds). We focus on the periodic orbits (oscillations) emerging from the Hopf bifurcation at a high value of ${A}$ (${A}\approx 15$, label (HB)). The black curves show the maximum and minimum values of $Y_3$ along this branch of mostly stable periodic orbits that emanate from this Hopf bifurcation and terminate at the low-${A}$ end in a homoclinic bifurcation of saddle-node on invariant circle type (SNIC \citep{izhikevich2007dynamical}), at which the frequency of oscillation goes to zero as ${A}$ approaches the value for the saddle node bifurcation on the lower branch of equilibrium points (${A}\approx 7$).

Figure~\ref{dimandnondim}b shows the frequency of these oscillations over the range of parameters $A$ where they exist (${A}\in[7,14.4]$). We observe a sharp transition between high-frequency oscillations for $A$ in the range $[10.2,14.4]$ where the frequency (in Hz) is in the range $[8,11]$, and low-frequency oscillation for ${A}$ in the range $[7,10.2]$, where the frequency is about $4$\,Hz and below. 
The oscillation frequencies on either side of the transition are associated with wakefulness state (alpha-like activity, around $10$\,Hz) and deep sleep (delta-like activity, around $2$\,Hz). 
%Jansen Rit \cite{jansen1995electroencephalogram} enables the model to generate  alpha similar to EEG alpha activity when they consider external input ($p(t)$, excitation  signal coming from thalamus or other distant column) as a random input uniformly distributed between 120 and 320 pulses per second. An initial occurrence of delta-like activity was reported in \cite{ahmadizadeh2018bifurcation}, where the authors provide detailed bifurcation analysis of two doubled cortical columns Jansen Rit model. 

Figures~\ref{dimandnondim}d and \ref{dimandnondim}e show that this change in frequency is accompanied with a qualitative change in the time profiles of the orbits. The alpha-type fast oscillations have a four-phase profile 
%\hmq{I do not understand?}
typical for an oscillatory negative feedback loop (compare Fig.~\ref{dimandnondim}d, and note that the scale differs between different outputs such that Fig.~\ref{dimandnondim}d has two $y$-axes):%\hmq{is this for alpha or delta? it is not clear }
\begin{itemize}
\item ($0$--$0.02$\,s) high pyramidal cell output ($Y_1$, blue) causes rising inhibition ($Y_2$, red),
\item ($0.02$--$0.05$\,s) high inhibition causes decrease in pyramidal cell output,
\item ($0.03-0.07$\,s) low pyramidal cell output  causes decreasing inhibition,
\item ($0.06-0.1$\,s) low inhibition causes rising pyramidal cell output.
\end{itemize} 
Throughout the entire period of the alpha-type oscillation the pyramidal cell output is supported by a near-constant input from the excitatory interneurons ($Y_3$, green). 
The delta-type (slow) oscillations are shown in Fig.~\ref{dimandnondim}e. They show a relaxation-type time profile: pyramidal cell output $Y_1$ (blue) and excitatory interneuron output $Y_3$ (green) collapse to zero and stay there for some period. This period is determined by how long it takes for inhibition ($Y_2$, red) to reach sufficiently low levels to permit pyramidal cell output $Y_1$ and excitatory interneuron output $Y_3$ to rise again.
The two time profiles in Fig.~\ref{dimandnondim}d and \ref{dimandnondim}e occur for parameter values of $A$ just above (Figure~\ref{dimandnondim}d, $A=11$) and below (Figure~\ref{dimandnondim}e, $A=10$) the transition. The bifurcation diagram shows even a tiny region of bistability, bounded by two saddle-node-of-limit-cycle bifurcations. However, this bistability and the saddle-nodes are not a consistent feature for this transition. Figure~\ref{dimandnondim}b shows a bifurcation diagram for non-zero external input ($p=120$\,Hz, yellow curve), where the transition still occurs and is still sharp, occurring in a small parameter region, but without bistability or saddle-nodes of limit cycles.

\citet{forrester2020role} detected this transition in their investigation of the Jansen-Rit model \eqref{jansenrit}. They demonstrated that it is an essential ingredient in the occurrence of large-scale oscillations in a network of neural populations (which would be measurable by EEG). As the transition is not associated with a bifurcation in parts of the parameter space (see Fig.~3 of \citep{forrester2020role}), %\hmq{they explore with zero input?} 
the authors labelled the transition a ``false bifurcation'' and tracked it in parameter space by associating it with a feature of the time profile, namely the occurrence of an inflection point.
\citet{forrester2020role} pointed out that ``false bifurcations'' are usually originating from singular perturbation effects in the system, referring to \citet{desroches2013inflection} and \citet{rodrigues2010method}.

To investigate this phenomenon further we non-dimensionalise the Jansen-Rit model \eqref{jansenrit} and identify a small parameter $\epsilon$ for which we can then study the singular limit for the transition from alpha to delta activity. Bifurcation diagram~\ref{dimandnondim}a shows a sharp change in the time profile of $Y_3$ from nearly constant for alpha activity to an oscillation with excursions close to zero: observe the drop in the minimum of the periodic orbit (black curve) for $A$ slightly above $10$.

The numerical bifurcation analysis  of the Jansen-Rit model is carried out  using the \textsc{coco} toolbox \citep{dankowicz2013recipes}, which is a MATLAB-based platform for parameter continuation allowing for bifurcation analysis of equilibria and periodic orbits. 
%We also compared all results of bifurcation analysis to simulations with the software \texttt{XXPAUT} \cite{ermentrout2003simulating}. 
Numerical continuation in \textsc{coco} \citep{dankowicz2013recipes} and \textsc{xppaut} \citep{ermentrout2002guide} is  used  to track stability of equilibria and periodic orbits and  to  detect their bifurcations. For numerical integration (simulation), we use  \textsc{xppaut} and \textsc{ode45} (a Runge–Kutta method)
in \textsc{matlab}. In all single-parameter bifurcation diagrams solid curves indicate stable states, while dashed lines are unstable states. Black curves indicate maximum and minimum values of periodic solutions in all  bifurcation diagrams. %  These properties are consistently maintained.
%We apply nondimensionalisation techniques in order to identify small parameter associated with  irregular behaviour of periodic orbit, as described in section \ref{}. This technique  allow us  to collect the constants in the problem so that the as many coefficients on the right hand side of the dimensionless equations  are of order $1$ as possible. This simplification allows us to more easily analyze the dynamics of the excitatory and inhibitory populations. By making these adjustments, we can obtain clearer insights into the behavior of the system. 

\section{Singular perturbation analysis of Jansen-Rit model  }
\label{sec:singpert:JR}
We introduce small parameter $\epsilon>0$ for which in the limit $\epsilon\to0$ the right-hand side (vector field) of system \eqref{jansenrit} becomes a piecewise linear non-smooth system.  First we non-dimensionalise the Jansen-Rit model~\eqref{jansenrit} to identify the small parameter.  
\subsection{Non-dimensionalisation of the Jansen-Rit model}\label{sec:nondimensionalisation}
 \citet{touboul2009codimension} present a non-dimensionalised version of the Jansen-Rit model. We choose the following dimensionless scale, as proposed by \citet{touboul2009codimension}:
\begin{itemize}
\item Dimensionless time %for excitatory population 
is set according to the internal decay time scale of the excitatory populations,
$t_\mathrm{new}=at_\mathrm{old}$, where ${a}$ is given in Table~\ref{tab:paramter model}.
\item We introduce the decay rate ratio ${b^*}$ of inhibitory to excitatory populations and  the ratio $G$ of postsynaptic amplitudes of inhibitory to excitatory populations 
\begin{align*}
  {b^*}&=\frac{{b}}{{a}},&G&=\frac{{B}}{{A}}.
\end{align*}
\item We rescale the state variables $Y_1$, $Y_2$, $Y_3$ such that they are all of order $1$ at their maximum, introducing a small parameter,
\begin{align}\label{def:epsilon}
 \epsilon&= \frac{2{a}}{{Br C(2e_0)}} \approx0.024\\
 y_1(t_\mathrm{new})&=rC\epsilon Y_1\left(\frac{t_\mathrm{new}}{a}\right),\ y_i(t_\mathrm{new})= r\epsilon Y_i\left(\frac{t_\mathrm{new}}{a}\right)\mbox{\quad for $i=2,3$.}\label{def:yscale}
\end{align}
\end{itemize}
We note that in Figure~\ref{dimandnondim}d and \ref{dimandnondim}e the quantities $Y_1$ and $Y_2$, $Y_3$ are on different $y$-axes, because they have different orders of magnitude. The non-dimensionalisation \eqref{def:yscale} ensures that pyramidal cell output $y_1$ has the same (order-$1$) magnitude as the inhibitory interneurons ($y_2$) and the excitatory interneurons ($y_3$).
The quantity $b^*$ is a measure of the difference in internal times scales between inhibitory and excitatory populations. Usually $b^*<1$ as inhibition is slower. The quantity $G$ is a measure for how strong internal feedback strength from inhibitory populations is compared to excitatory populations. 

%{\centering $y_i(\tau)= r\epsilon{Y_i}(\frac{\tau}{a})$, with $i=2,3$  for  PSPs of the interneuron populations.}

By substituting  these dimensionless dependent variables   into  equation \eqref{jansenrit}, and applying the chain rule with respect to the dimensionless time scale $t_\mathrm{new}$, we obtain the  dimensionless form (using $\Dot{y}_i$ and $\Ddot{y}_i$ also for the new time)

\begin{align} 
 \begin{split} 
\ddot y_1&=\frac{2}{G}\newS_{\epsilon/a_1}\left(y_3-y_2-y_{0,1}\right)-2\dot y_1-y_1, \\ 
\ddot y_2&=2b^*\alpha_4\newS_{\epsilon/a_2}(y_1-y_{0,2})-2b^{*} \dot y_2-(b^{*})^2y_2,\\
\ddot y_3&=\frac{P}{G}+\frac{2\alpha_2}{G}\newS_{\epsilon/a_3}(y_1-y_{0,3})-2 \dot y_3-y_3,\\
\end{split}
\label{nondimwitheps}
\end{align}
where new dimensionless sigmoidal transformation is (see Fig.~\ref{fig:sigmwithdiffeps})
\begin{align}  
{\newS_{\epsilon}
(y)}=\cfrac{1}{1+\exp \left(-y/\epsilon\right)}.
\label{newsigwitheps}
\end{align}
\begin{table}[ht]
    \centering
     \caption{Parameter values of the dimensionless model equation \eqref{nondimwitheps}.} 
    \begin{tabular}{lP{0.4\linewidth}P{0.33\linewidth}}
        \toprule % Thicker line at the top
        Parameter & Relation to the original parameters & Value \\ 
        \midrule % Line between header and content
%        $\beta$  & Scaling factor   & 0.025   \\
        $\epsilon$  & inverse of  slope of activation function for excitatory populations %\hmq{1/4e}  
        & 0.024   \\
        $y_{0,1},y_{0,2},y_{0,3}$   & re-scaled activation thresholds   &{$y_{0,1}=y_{0,3}=0.08$,} {$y_{0,2}=0.3$}\\
       $G$ &$B/A$ ratio of feedback strengths &$[0.5,20]$ (varied)\\
        $b^*$ & $b/a$ time scale ratio (inh. vs. exc.)&[0.2,0.5] (varied, original $0.5$)\\
%$jnew$& $\frac{Br\beta C(2e_0)}{a}$&2.079\\
$P$& $\epsilon Br p(t)/a$, external input&0\\
$\alpha_2,\alpha_4$& average number of synaptic connections (same as original)&$\alpha_2=0.8$, $\alpha_4=0.25$ \\
$a_i$& population-dependent factors in sigmoid slope& $a_1\!=\!a_3\!=\!1$, $a_2\!=\!1/4$.\\
        \bottomrule % Thicker line at the bottom
    \end{tabular}
    \label{tab:nondimparameters}
\end{table}

As we have rescaled the PSP's $Y_i(t)$, the thresholds in the activation function $\newS$ are now different for each neuron population, such that we call them $y_{0,i}$. These thresholds $y_{0,i}$, at which the activation $\newS_{\epsilon/a_i}(y-y_{0,i})$ equals $1/2$, now show up in \eqref{nondimwitheps}.
The new non-dimensional activation thresholds are  \begin{align}\label{def:nondim:y0}
    y_{0,1}&=\frac{ry_0}{a_1}\epsilon,& y_{0,2}&=\frac{ry_0}{a_2}\epsilon,& y_{0,3}&=\frac{ry_0}{a_3}\epsilon, \mbox{\quad where $a_1=a_3=1$, $a_2=1/4$,} 
\end{align}
using the parameters from Table~\ref{tab:paramter model}, which result in the non-dimensional parameters shown in Table~\ref{tab:nondimparameters}. 
Note that the main difference to Touboul's non-dimensionalisation \citep{touboul2009codimension} is that the definition of the sigmoid function includes the scaling factor $\epsilon$.
Figure~\ref{dimandnondim}c shows the same  bifurcation diagram as shown in Figure~\ref{dimandnondim}a, but uses  the non-dimensionalized quantities $G$ and $y_3$.
The new primary bifurcation parameter $G$ (the ratio of feedback strengths between inhibition and excitation) is proportional to the inverse of the excitation feedback strength $A$, such that now the alpha-to-delta frequency transition occurs for increasing $G$ at $G\approx2.2$. The Hopf bifurcation occurs for low $G\approx1.5$ and the SNIC connecting orbit occurs at $G\approx3$.

\subsection{\texorpdfstring{Discussion of smallness of system parameters}{Discussion of smallness of system parameters}}\label{sec:whatissmall}
In the non-dimensionalised model \eqref{newsigwitheps} the small parameter $\epsilon$, which equals $0.024$, appears in the
inverse of the slope of the dimensionless sigmoid $\newS_\epsilon$ at $y=0$ %\hmq{? the slop is 1/4epsilon}, 
letting the sigmoid $\newS_\epsilon$ approach a discontinuous switch in the limit for small $\epsilon$ and $i=1,2,3$:
\begin{align} 
 \lim_{\epsilon \to 0} \newS_{\epsilon}
(y-y_{0,i}) &=
 \begin{cases}
		0  & \mbox{if  $y<y_{0,i} $,} \\
		1 & \mbox{if  $y>y_{0,i}$,}
\end{cases}\label{switches}
\end{align}
In addition the activation thresholds $y_{0,i}$ inherit a factor $\epsilon$ in \eqref{def:nondim:y0}. However, we observe that there are non-trivial factors in front of $\epsilon$ in several places. Since the $\epsilon$ in the original parameters setting is about $1/40$ (so not that small), these factors will influence our strategies for taking the singular limit $\epsilon\to0$.
\begin{itemize}
\item The slope in the activation for the inhibitory interneurons, $S_{\epsilon/a_2}$, is$a_2/(4\epsilon)=1/(16\epsilon)\approx 2.6$. So, it is further away from the limiting discontinuous switch than the activation of the excitatory populations.
\item
The different factors $ry_0/a_1$, $ry_0/a_3$ and $ry_0/a_2$ of $\epsilon$ in \eqref{def:nondim:y0} are obtained by  substituting the baseline values of the original parameters from Table \ref{tab:paramter model} into system \eqref{nondimwitheps}. We observe that the activation thresholds of excitation are small ($y_{0,1}=y_{0,3}\approx 0.08$), but the threshold of inhibition  $y_{0,2}$ ($y_{0,2}\approx0.3$) is not a small quantity.
%\item The factors $ry_0/a_1$, $ry_0/a_3$  are such that $y_{0,1}=y_{0,3}\approx 0.08$, such that the activation thresholds are small, but not as small as $\epsilon$.
%\item The factor $ry_0/a_2$ is such that $y_{0,2}\approx 0.3$, such that $y_{0,2}$ is not a small quantity.
\end{itemize}
The smallness of $\epsilon$ expresses that the internal dynamics of neurons in the excitatory populations is fast, leading to small thresholds and steep activation functions after non-dimensionalisation. The inhibitory neural population is comparatively slower, leading at the population level to a shallower slope of the activation curve and a larger activation threshold  \citep{wilson1972excitatory}.
The above observations suggest  several possible singular limits.
\begin{itemize}
    \item[(a)] We let $\epsilon$ go to zero in the denominator appearing in $\newS_{\epsilon/a_i}$ for all neuron populations ($i=1,2,3$) simultaneously. At the same time we keep the activation thresholds $y_{0,i}$ fixed for all populations. This leaves the parameters $y_{0,i}$ in the model as independent parameters and uses the concrete values from Table~\ref{tab:nondimparameters}.
    \item[(b)] We let $\epsilon$ go to zero in the denominator appearing in $\newS_{\epsilon/a_i}$ for all neuron populations ($i=1,2,3$) simultaneously, and we let the activation thresholds $y_{0,1}$ and $y_{0,3}$ for the excitatory neuron populations go to zero (either proportional to $\epsilon$ or at some lower rate).
    \item[(c)] We let $\epsilon$ go to zero for the excitatory neurons populations ($i=1,3$) in the slope of the activations $\newS_{\epsilon/a_i}$, but keep the activation for the inhibition with a fixed finite slope, equal to $\approx 10$ (as well as keeping the activation threshold $y_{0,2}\approx0.3$ fixed). 
\end{itemize}
We will focus in our analysis on strategy (a) in this paper, because it permits explicit expressions for most equilibria and their bifurcations and it is sufficient to derive in an implicit algebraic condition for the  alpha-to-delta transition. In both limits the collapse of the alpha-frequency oscillations is a \emph{grazing bifurcation} \citep{bernardo2008piecewise} of periodic orbits in a piecewise linear ODE. In limit (a) the collapse leads to a low excitation equilibrium ($y_i\ll1$ for $i=1,2,3$), in limit (b) it leads to delta activity. 

In the numerical bifurcation diagram in Fig.~\ref{dimandnondim}c the three boundaries for alpha and delta activity are the Hopf bifurcation (onset of alpha frequency oscillations) at $G\approx 1.5$, the transition between alpha and delta activity at $G\approx 2.2$, and the connecting orbit to a saddle-node (SNIC bifurcation) at $G\approx 3$. Two of the boundaries can be found by analysis of equilibria, namely the Hopf bifurcation and the saddle-node bifurcation for small $y_1$.

For our analysis, the dimensionless external input $p$ is fixed to zero.   
We use the dimensionless parameter $G$ (the ratio between feedback strength of inhibition versus excitation) as our primary bifurcation parameter. We will later vary $b^*$, the ratio between internal time scales between inhibition and excitation as a secondary bifurcation parameter.

\subsection{Equilibrium analysis  in the small-\texorpdfstring{$\epsilon$}{epsilon} limit}
\label{subsec:Equilibria in the  limit}

Inserting $\epsilon=0$ (or $\epsilon\ll1$) into the activation functions in \eqref{nondimwitheps} and \eqref{newsigwitheps}, results in an easy-to-interpret limit for equilibria. 
Setting all derivatives in equation \eqref{nondimwitheps} to zero, we obtain that the equilibrium values for $y_1$, $y_2$, $y_3$ satisfy the algebraic equations
\begin{align} 
 y_1&=\frac{2}{G}\newS_\epsilon(y_3-y_2- y_{0,1})\label{xeqy1}\\   
y_2&=\frac{2\alpha_4} {b^*}\newS_{\epsilon/a_2}(y_1-y_{0,2}) \label{xeqy2}\\y_3&=\frac{2\alpha_2}{G}\newS_\epsilon(y_1-y_{0,3}).\label{xeqy3}
\end{align}
The right-hand sides of \eqref{xeqy2} and \eqref{xeqy3} define the equilibrium values $y_2$ and $y_3$ as a function of $y_1$. This is also true in the limit $\epsilon=0$, whenever $y_1$ is not on the activation threshold ($y_{0,2}$ and $y_{0,3}$, respectively).

Inserting \eqref{xeqy2} and \eqref{xeqy3} into \eqref{xeqy1} and taking the Heaviside limit $S_0$ in \eqref{switches} for the  activation functions $S_{\epsilon/a_i}$
(indicating the dependence of $y_2$ and $y_3$ on $y_1$), leads to the relation
\begin{align}
   \frac{G}{2}y_1 =\rhs_0(y_1):= \begin{cases}
		0  & \mbox{if  $y_3(y_1-y_{03})-y_2(y_1-y_{02})< y_{01}$,} \\
		1& \mbox{if  $ y_3(y_1-y_{03})-y_2(y_1-y_{02})> y_{01}$,}
\end{cases} \label{y1piecewise}
\end{align}
\begin{figure}[ht] 
    \centering
\includegraphics[width=\textwidth
]{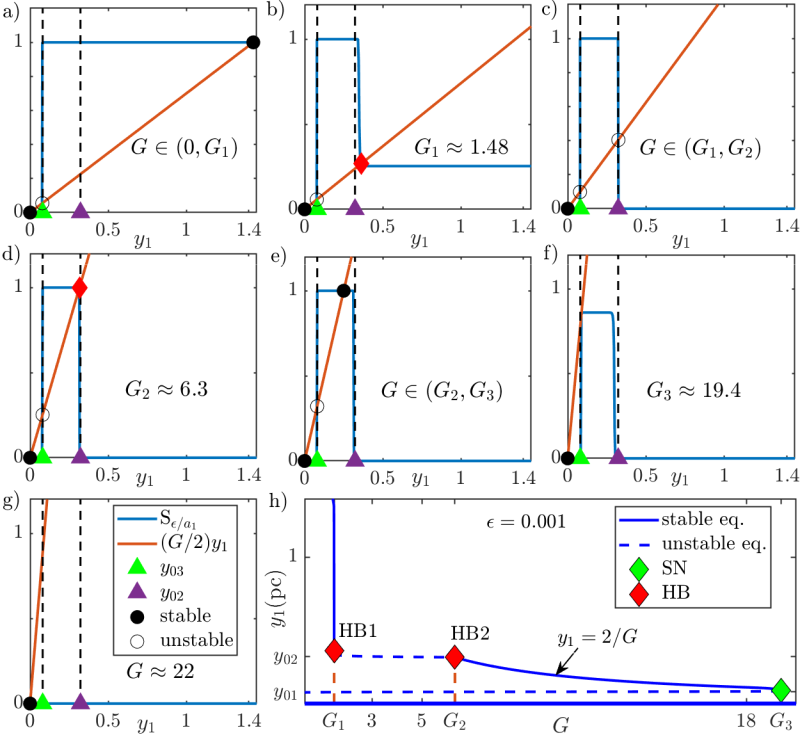}
    \caption{(a-g) Intersections between left- (red) and right-hand (blue) side of \eqref{y1piecewise} for different bifurcation parameters $G$ (plots used $\epsilon=0.001$, $y_{0,1},y_{0,3}=0.08$  and $y_{0,2}=0.3$). (h) Numerical bifurcation diagram of equilibria in $(G,y_1)$-plane obtained by  \textsc{coco} \citep{dankowicz2013recipes}. $\newS_\epsilon/a_1=\newS_\epsilon(y_3-y_2- y_{0,1})$ for $y_1$ in equation \eqref{xeqy1}.
    For other parameters see Table~\ref{tab:nondimparameters}.} 
    \label{intersectionfig}
\end{figure}
\begin{table}[ht]
\centering
\caption{\label{tab:orgionalthresholds}Equilibrium points present in system \eqref{nondimwitheps} for all parameter ranges of $G$: for critical values of $G$ see \eqref{def:g1g3}, \eqref{def:g2g4}. The index is the number of unstable eigenvalues of the Jacobian in the equilibrium for non-zero $\epsilon$.}
       \begin{tabular}{lcP{0.35\linewidth}c}
        \toprule % Thicker line at the top
     range of $G$&Fig.\,\ref{intersectionfig}&equilibrium location & index \\
      \midrule % Line between header and content
         $I_1=(0,G_1)$  &(a) & $y_\mathrm{eq_{1}}=(0,0,0)$ & $0$\\
         &&$y_\mathrm{eq_{2}}=(y_{0,3},0,y_{0,1})$& $1$\\
         &&$y_\mathrm{eq_{3}}=(2/G,2\alpha_4/b^*,2\alpha_2/G)$&$0$ \\
    $G_1$: Hopf bifurcation& (b)\\
    $I_2=(G_1,G_2)$ & (c) &$y_\mathrm{eq_{1}}=(0,0,0)$ & $0$\\
    &&$y_\mathrm{eq_{2}}=(y_{0,3},0,y_{0,1})$& $1$\\
    &&$y_\mathrm{eq_{4}}=(y_{0,2},2\alpha_2/G-y_{0,1},2\alpha_2/G)$&$2$\\
    $G_2$: Hopf bifurcation& (d)\\    
    $I_3=(G_2,\min(G_3,G_4))$&(e)& $y_\mathrm{eq_{1}}=(0,0,0)$& $0$\\
    && $y_\mathrm{eq_{2}}=(y_{0,3},0,y_{0,1})$&$1$\\
    && $y_\mathrm{eq_{4}}=(2/G,0,2\alpha_2/G)$ &$0$\\
    $\min(G_3,G_4)$: saddle-node&(f)\\   
    $I_4=(\min(G_3,G_4),\infty)$&
    (g)& $y_\mathrm{eq_{1}}=(0,0,0)$ &$0$ \\       
  \bottomrule % Thicker line at the bottom
  \end{tabular} 
  \label{tab:xeqresults}
\end{table}

Analysis of \eqref{xeqy1}--\eqref{xeqy3} in the limit $\epsilon\to0$ is supported by Figure~\ref{intersectionfig}. Each panel shows the left-hand side of \eqref{y1piecewise} (in red) and the right-hand side $\rhs_0(y_1)$ of \eqref{y1piecewise} (in blue) as a function of $y_1$ for different values of $G$, otherwise using the parameter values in Table~\ref{tab:nondimparameters}. The (red) left-side function is the straight line $y_1\mapsto Gy_1/2$. The right-side function $\rhs_0(y_1)$ is piecewise constant. Defining two critical values for parameter $G$ as
\begin{align}
    \label{def:g1g3}
    G_1&=\frac{\alpha_2b^*}{\alpha_4+y_{0,1}b^*/2}\approx1.48,&
    G_3&=\frac{2\alpha_2}{y_{0,1}}\approx19.4,
\end{align}
(thus, $0<G_1<G_3$) the right-hand side equals (recall that $0<y_{0,3}<y_{0,2}$)
\begin{equation*}
\begin{aligned}
    \\
    \rhs_0(y_1)&=
    \begin{cases}
    \\
    \\
    \
    \end{cases}
\end{aligned}\hspace*{-2em}
\begin{aligned}
    &&\mbox{for $G$ in:}&&(0,&G_1),&
    (G_1,&G_3),&(G_3,&\infty),\\
    y_1&\in(0,y_{0,3})& \mbox{($y_2$ off, $y_3$ off):}&& &0 &&0 &0\\
    y_1&\in(y_{0,3},y_{0,2})& \mbox{($y_2$ off, $y_3$ on):}
    && &1 &&1 &0\\
    y_1&\in(y_{0,2},\infty)&
    \mbox{($y_2$ on, $y_3$ on):}
    && &1 &&0 &0\\
    \end{aligned}
\end{equation*}
%\hmq{whats is 1,0 ?}
Each intersection of the two lines in Fig.~\ref{intersectionfig} is a stable equilibrium if $\rhs_0$ (blue line) is horizontal in the intersection, or an unstable, so-called \emph{pseudo-equilibrium} \citep{kuznetsov2003one,bernardo2008piecewise}, if $\rhs_0$ is vertical in the intersection. For $\epsilon\ll1$ but non-zero, the pseudo-equilibria of the piecewise linear limit system will become unstable equilibria with strongly unstable eigenspaces.
Table~\ref{tab:xeqresults} lists the locations of equilibria for all ranges of parameter $G$, with two additional critical values for $G$,
\begin{align}
\label{def:g2g4}
    G_2&=\frac{2}{y_{0,2}}\approx6.3,&
    G_4&=\frac{2}{y_{0,3}}\approx24.8.
\end{align} Figure~\ref{intersectionfig}h shows the resulting bifurcation diagram in the limit of small $\epsilon$.

Figure~\ref{intersectionfig}a shows that for $G<G_1$ we have a pair of a stable and unstable equilibria, both with small values of $y_1$: the unstable equilibrium is a saddle with $y_1=y_{0,3}=0.08$ for $0<\epsilon\ll1$, and the stable equilibrium is a node with a $y_1$-value very close to zero, on the order $\exp(-y_{0,1}/\epsilon)$.

%Figure~\ref{intersectionfig}a shows that for $G<G_1$ we have a pair of a stable and unstable equilibria with $y_1$ near $0$, where the saddle is has $y_1=y_{0,3}$, and the node is very close to zero $\sim \exp(-y_{0,1}/\epsilon)$). 
There exists another stable equilibrium, for which  $y_1$ is large, equal to $2/G$. At $G=G_1$ (Fig.~\ref{intersectionfig}b, see \eqref{def:g1g3} for value of $G_1$) $\rhs_0$ drops to zero for $y>y_{0,2}$, which results in a sudden drop of the $y_1$-component of the large-$y_1$ equilibrium to $y_{0,2}$ and its loss of stability. For positive $\epsilon$ this loss of stability is a Hopf bifurcation (see the labeled point HB1 in Fig. 3h). Figure~\ref{intersectionfig}c shows the equilibria for $G\in(G_1,G_2)$. Figure~\ref{intersectionfig}d shows the situation at the next critical  value of $G$, the Hopf bifurcation $G=G_2=2/y_{0,2}$,  when the left-hand side function (red line) goes through the ``corner'' of $\rhs_0$ at $y_1=y_{0,2}$ (see red diamond). At this value the equilibrium with $y_1=y_{0,2}$ regains its stability in a Hopf  bifurcation for non-zero $\epsilon$ (see the labeled point HB2 in Fig. 3h). Fig.~\ref{intersectionfig}e shows the situation for $G>G_2$ with two stable equilibria and one saddle. The next critical $G$ is $\min(G_3,G_4)$, when either the excitatory interneuron activity is no longer sufficient to overcome the threshold ($G_3$, when the blue line in Fig.~\ref{intersectionfig}g  drops to zero on the $y_1$-interval $(y_{0,3},y_{0,2})$), or  the left-hand side (red curve) touches the ``corner'' of $\rhs_0$ at $y_1=y_{0,3}$ ($G_4$). In both cases, the resulting equilibrium bifurcation is a saddle-node bifurcation of the equilibrium with$y_1=y_{0,1}=y_{0,3}$ (see the labeled SN  in Fig.~\ref{intersectionfig}h). For the parameters we study, we have that $G_3<G_4$ such that $G_4$ does not change the dynamics. % \cor{check with Jan as there is no SN (saddle xeq?)}

\paragraph{Brief comment on the excitatory activation thresholds $y_{0,1}$ and $y_{0,3}$} 
While we will treat the excitatory thresholds $y_{0,1}$ and $y_{0,3}$ as positive constants in our analysis in Section~\ref{sec:Piece-wise periodic orbits} (called strategy (a) in Section~\ref{sec:whatissmall}), let us perform a quick check what happens if we also assume that the excitatory thresholds $y_{0,1}$ and $y_{0,3}$ are assumed to be small. If we assume proportionally small excitatory activation thresholds $y_{01}=y_{03}\sim \epsilon$, then the saddle and node equilibria with $y_1\approx0$ are no longer present for any $G$ of order $1$. This follows immediately from a perturbation analysis of the equilibria for small non-zero $\epsilon$ and excitatory activation thresholds of the form 
\begin{align*}
 y_{0,1}&=\epsilon z_{0,1}, &  
 y_{0,3}&=\epsilon z_{0,3}, &  
 z_{0,1},z_{0,3}&>0, &
 z_{0,1},z_{0,3}&=O(1).
\end{align*}
Let us assume that $y_{1,\mathrm{eq}}\ll1$, in particular, $y_{1,\mathrm{eq}}\ll y_{0,2}$, and check that no such equilibrium can exist. The smallness of $y_{1,\mathrm{eq}}$ together with identity \eqref{xeqy2}, $y_{2,\mathrm{eq}}=2\alpha_4/(b^*(1+\exp((y_{0,2}-y_{1,\mathrm{eq}})/(\epsilon/a_2)))$, would imply that the corresponding inhibitory interneuron activity $y_{2,\mathrm{eq}}$ is very close to zero: $y_{2,\mathrm{eq}}\sim \exp(-a_2 y_{0,2}/\epsilon)\leq\epsilon$. From this, the identities \eqref{xeqy1} and \eqref{xeqy3} imply positive  lower bounds independent of $\epsilon$ for the equilibrium pyramidal-cell activity $y_{1,\mathrm{eq}}$ and the excitatory interneuron activity $y_{3,\mathrm{eq}}$, namely 
\begin{align} \label{smallthresh:ylowbd}
\begin{split}
    y_{1,\mathrm{eq}}&\geq\frac{2}{G}\frac{1}{1+\exp(z_{0,1}+y_{2,\mathrm{eq}}/\epsilon)}\geq \frac{2}{G}\frac{1}{1+\exp(z_{0,1}+1)},\\
    y_{3,\mathrm{eq}}&\geq\frac{2\alpha_2}{G}\frac{1}{1+\exp(z_{0,3})}
\end{split}
\end{align}
Hence, for excitatory activation thresholds of order $\epsilon$ the pair of small-$y_1$ equilibria does not exist. 
The other equilibria in Table~\ref{tab:xeqresults} and Fig.~\ref{intersectionfig}h have well-defined limits for $y_{0,1}, y_{0,3}\to0$: the stable equilibrium has the limit $y_1=2/G$, the unstable equilibrium has the limit $y_1=y_{0,2}$ and its inhibitory activity  $y_2$ approaches $2\alpha_2/G$ for small $y_{0,1}$. The saddle-node bifurcation at $\min(G_3,G_4)$ goes to infinity as $y_{0,1},y_{0,3}\to0$, such that the large-activity equilibrium exists over a wider range of $G$. This leaves the region $(G_1,G_2)$
between the two Hopf bifurcations, where no stable equilibrium exists for $y_{0,1},y_{0,3}$ of order $\epsilon$. The Hopf bifurcation at $G_1$ has the limit $\alpha_2b^*/\alpha_4$ for $y_{0,1}\to0$, while $G_2$ is independent of $y_{0,1}$ and $y_{0,3}$.
\subsection{Periodic orbits in the small-\texorpdfstring{$\epsilon$}{epsilon} limit}
\label{sec:Piece-wise periodic orbits}

\paragraph{The case $\epsilon=0$ as a piecewise smooth ODE} System \eqref{nondimwitheps} in the limit $\epsilon\to0$ has discontinuities on the right-hand side such that it is only piecewise smooth in its phase space \citep{bernardo2008piecewise}. For such systems there is in general no guarantee that trajectories can be consistently continued across discontinuities by following the flow defined by the right-hand sides in each part of the phase space. The problem can be illustrated for simplest scenario, an ODE with a right-hand side consisting of two smooth pieces and a single discontinuity, 
\begin{align*}
    \dot x(t)=\begin{cases}
        f_+(x)&\mbox{if $h(x)\geq0$,}\\
        f_-(x)&\mbox{if $h(x)<0$}
    \end{cases}
    \end{align*} 
    with smooth right-hand sides $f_\pm:\R^n\to\R^n$ and \emph{switching function} $h:\R^n\to\R$ (e.g.,  consider the scalar example $\dot x(t)=-\sign x(t)$, where $f_\pm(x)=\mp1$ and $h(x)=x$). For all points with $\pm h(x)>0$ it is clear that the trajectory through $x$ is determined by $f_\pm(x)$. However, when following a trajectory $x(t)$ one may approach a point $x$ with $h(x)=0$ (on the \emph{discontinuity surface} $\{x:h(x)=0\}$). If in this point $x$ the quantities $h'(x)f_-(x)$ and $ h'(x)f_+(x)$ have opposite sign (using $h'(x)$ for $\partial h(x)$), then the trajectory cannot ``cross to the other side'' of the discontinuity surface $D=\{x:h(x)=0\}$ consistent with the right-hand sides. A commonly adopted convention for this case where the flows on both sides of a codimension-$1$ discontinuity surface point in opposite directions relative to the surface, is continuing the trajectory using the \emph{Filippov} solution, or \emph{sliding}, by following the flow of \begin{align*}
    \dot x(t)&=\frac{f_+(x)h'(x)f_+(x)-f_-(x)h'(x)f_-(x)}{h'(x)f_+(x)-h'(x)f_-(x)}&
    \mbox{if\quad} \begin{cases}
        \mbox{$h(x)=0$ and}\\
            \mbox{$h'(x)f_+(x)h'(x)f_-(x)<0$.}
    \end{cases}
\end{align*} 
However, the structure of \eqref{nondimwitheps} rules out this common scenario of sliding, where $h'(x)f_+(x)h'(x)f_-(x)<0$, for the discontinuity surfaces in system \eqref{nondimwitheps} generated by the limit $\epsilon\to0$ in \eqref{newsigwitheps} System \eqref{nondimwitheps} has $3$ hyperplanes where it is discontinuous in the limit $\epsilon\to0$, with the switching functions
\begin{align*}
    h_1(y,\dot y)&=y_y-y_2-y_{0,1},&h_2(y,\dot y)&=y_1-y_{0,2},&h_3(y,\dot y)&=y_1-y_{0,3}.
\end{align*}
In each of the hyperplanes $\{h_i=0\}$, where the vector field is discontinuous, the two half-space vector fields point into the same direction relative to the hyperplane. This is due to the presence of inertia: \eqref{nondimwitheps} is a system of second-order equations for the firing activities but the switches in $\newS_0$ depend only on the firing activities, not their time derivatives. Let us inspect each codimension-$1$ discontinuity set:
\begin{itemize}
    \item $D_i=\{y:y_1=y_{0,i}\}$ for $i=2,3$, such that the switching function is $h_i(y,\dot y)=y_1-y_{0,i}$: the inner product of the vector fields $(\dot y,\ddot y)$ on both sides of $D_i$ with the gradient of $h_i$ equals $\Dot{y}_1$, which is continuous  across $D_i$ and, thus, all non-equilibrium trajectories cross the surface $D_i$ in the same direction from both sides of $D_i$.
    \item $D_1=\{y:y_3-y_2=y_{0,1}\}$, such that the switching function is $h_i(y,\dot y)=y_3-y_2-y_{0,1}$: the inner product of $(\dot y,\ddot y)$ on both sides of $D_1$ is $\Dot{y}_3-\Dot{y}_2$, which is also continuous across $D_1$.
\end{itemize}
Hence, the set $\mathcal{D}$ in phase space where trajectories are not well defined in the classical sense is much smaller than for common piecewise smooth ODEs. It is a subset of the codimension-$2$ surfaces $D_i'=\{(y,\dot y):y_1=y_{0,i},\dot y_1=0\}$ for $i=2,3$ and $D_1'=\{(y,\dot y):y_2-y_3=y_{0,1}, \dot y_2=\dot y_3\}$. Examples of points in these sets $\mathcal{D}$ are the ``unstable equilibria'' $y_{\mathrm{eq}_2}$ listed in Table~\ref{tab:xeqresults}. These types of equilibria are called pseudo equilibria in the literature about Filippov systems. As the equilibrium discussion in Figure~\ref{intersectionfig} shows, the pseudo equilibria for $\epsilon=0$ are limits of unstable equilibria for $\epsilon>0$.

\begin{figure}[ht]
    \centering
    \includegraphics[width=0.45\textwidth]{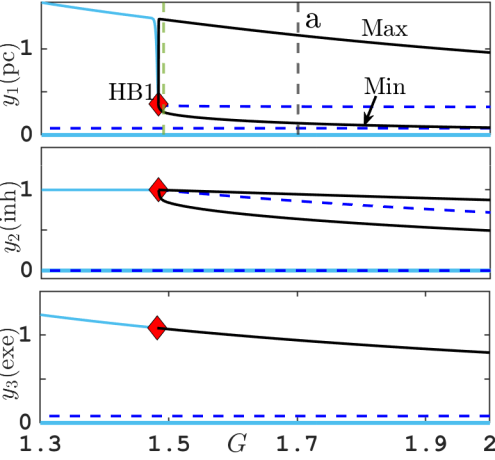}
    \caption{Equilibria and periodic orbits branching off Hopf bifurcation (HB1) for system~\eqref{nondimwitheps} for varying $G$ and small $\epsilon=0.001$, showing equilibrium values, or maxima and minima of the neural activities $y_1,y_2,y_3$, respectively. Green vertical line corresponds to canard orbit shown in  Fig~\ref{canard}. Grey  vertical line is at parameter $G=1.7$, used in singular limit in Fig.~\ref{figsamllepstime}(a). Other parameters are listed in Table \ref{tab:nondimparameters}.}
    \label{figbdsamlleps}
\end{figure}

\paragraph{Numerical bifurcation analysis of $\alpha$-type periodic orbits for small positive $\epsilon$}  As seen in Table~\ref{tab:xeqresults} and Figure~\ref{intersectionfig}h, in the small-$\epsilon$ limit  the high-activity equilibrium $y_1$ changes not only its stability at $G=G_1$ but also its location: the equilibrium pyramidal-cell activity drops sharply from $y_1=2/G=2/G_1\approx1.35$ to $y_1=y_{0,2}\approx0.3$. Hence, in the singular limit we cannot expect the standard Hopf bifurcation scenario with a family of nearly harmonic small-amplitude periodic orbits branching off. 

We initially perform a numerical continuation for the full non-dimen\-sional\-ised system \eqref{nondimwitheps} for small $\epsilon$ ($\epsilon=0.001$). Figure~\ref{figbdsamlleps} zooms into the range of parameters $G$ where periodic orbits of $\alpha$ type exist, near the Hopf bifurcation at $G=G_1\approx1.48$ (see \eqref{def:g1g3}). We observe that the pyramidal-cell ($y_1$) amplitude for the family of periodic orbits emerging from the Hopf bifurcation grows over a small range of $G$ to order $1$ (black curve in top panel in Figure~\ref{figbdsamlleps}), by approximately the same amount as the equilibrium value of $y_1$ dropped before the Hopf bifurcation. This explosive growth in amplitude is reminiscent of \emph{canard explosions} as first described for the classical van der Pol or FitzHugh-Nagumo oscillator by  \citet{benoit1981chasse} (see \citet{krupa2001relaxation} for general theory and \citet{wechselberger2007canards} for a review). However, the inhibition amplitude $y_2$ grows square-root-like with $G-G_1$, as occurs at a classical Hopf bifurcation. This discrepancy and the fact that the singular limit $\epsilon\to0$ is not slow-fast but discontinuous suggests that another mechanism must be responsible for the explosive growth in amplitude of $y_1$. Exploration of this phenomenon is beyond the scope of this paper. The time profile of periodic orbits near the explosion is shown in Figure~\ref{canard} in the discussion in Section~\ref{sec:discussion}.

\paragraph{Piecewise linear oscillations in a negative feedback loop between pyramidal cells and inhibition for $\epsilon=0$}
In Figure~\ref{figbdsamlleps} we also note that for $G\approx G_1$ and pyramidal-cell activity $y_1$ in the range between $y_{0,2}$ and $2/G_1$ the excitatory interneurons are uniformly active, as $y_1$ is always above the threshold $y_{0,3}$ at which excitatory interneuron activity is switched off, and which is much smaller than $y_{0,2}$. Hence, $y_3$ will be positive and at its equilibrium value  \begin{align*}
    y_3(t)\approx y_{3,\mathrm{eq}}=\frac{2\alpha_2}{G}\approx \frac{2\alpha_2}{G_1}= y_{0,1}+\frac{2\alpha_4}{b^*}.
\end{align*}
Thus, we may replace $y_3$ in the sigmoid input for $y_1$ with its equilibrium value, resulting in the four-dimensional system for pyramidal-cell and inhibitory activity only, inserting the limit $\epsilon=0$: 
\begin{align} 
     \mbox{(if\ }y_1(t)&>y_{0.3}\mbox{ for all time)}\nonumber\\
	\ddot y_{1}&= (2/G)\,u_1-2\dot y_{1}-y_{1},\label{y1dot}\\ 
 \ddot y_{2}&= 2b^*\alpha_4\ u_2-2b^*\dot y_{2}-(b^*)^2y_{2}\mbox{,\quad where}\label{y2dot}\\
     u_1(t)&=\begin{cases}
     1&\mbox{if $y_2<2\alpha_2/G-y_{0,1}$,}\\
     0&\mbox{if $y_2>2\alpha_2/G-y_{0,1}$,}
 \end{cases}\quad
 u_2(t)=\begin{cases}
     0&\mbox{if $y_1<y_{0,2}$,}\\
     1&\mbox{if $y_1>y_{0,2}$.}     
 \end{cases}
 \nonumber
 \end{align} 
System~\eqref{y1dot},\,\eqref{y2dot} is a classical negative feedback loop with inertia: $y_2$ suppresses $y_1$, while $y_1$ promotes $y_2$. Each oscillation period has four phases and follows a piecewise linear second-order vector field in each of the phases. 
\begin{figure}[ht]
    \centering
    \includegraphics[width=\textwidth]{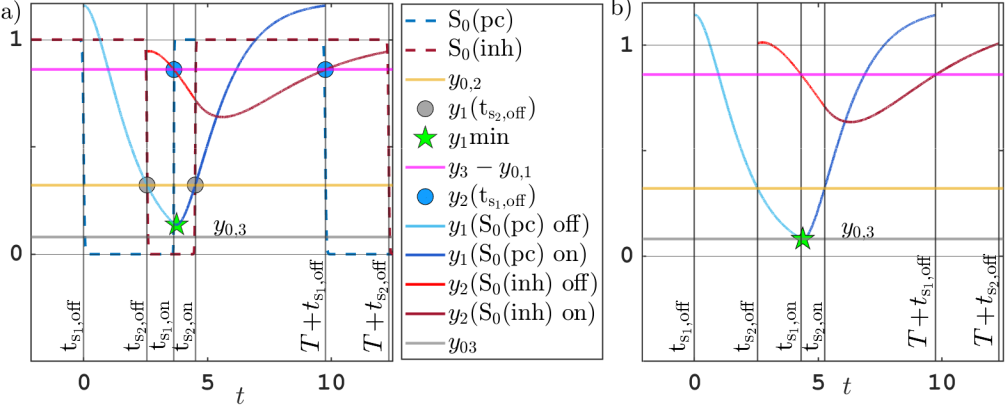}
     \caption{Time profile of alpha-type piecewise exponential periodic orbits given in \eqref{pworbit:crossing},representing solutions of the piecewise linear ODE~\eqref{y1dot},\,\eqref{y2dot}  for $y_{0,1},y_{0,3}=0.08$, $y_{0,2}=0.3$, and $y_3$ fixed at its equilibrium value $2\alpha_2/G$. Panel (a): orbit for  $(b^*,G)=(0.5,1.7)$; $t_\mathrm{s_2,off}$, $t_\mathrm{s_2,on}$: threshold crossing times of $y_1$ (where $y_1=y_{0,2}$); $t_\mathrm{s_1,off}$, $t_\mathrm{s_1,on}$: threshold crossing times of $y_2$ (where $y_2=y_3-y_{0,1}$). Dashed blue line (legend entry $S_0$(pc)) is activation switch $\newS_{0}\left(y_3-y_2(t)-y_{0,1}\right)$ for $y_1$ with $\epsilon=0.001$. Dashed red line (legend entry $S_0$(inh)) is activation switch $\newS_{0}\left(y_1(t)-y_{0,2}\right)$ for $y_2$. Panel (b): same as panel (a) but for $(b^*,G)=(0.44,1.7)$, where orbit grazes ($y_{1,\min}=y_{0,3}$, green star on grey horizontal line $\{y=y_{0,3}\}$).
    For other parameters see Table \ref{tab:nondimparameters}.}
    \label{figsamllepstime}
\end{figure}
Figure~\ref{figsamllepstime}a shows the typical shape for the $y_1$ and $y_2$ profile of the oscillations. We reduce the periodic problem for the four-dimensional piecewise affine ODE to a four-dimensional algebraic problem, parametrising the oscillations using the times at which the respective activations $u_1, u_2$ switch between $0$ and $1$ (w.l.o.g. $t_\mathrm{s_1,off}=0$), see Figure~\ref{figsamllepstime}a for guidance:
\begin{equation}\label{def:t_switch}
\begin{aligned}
    t_\mathrm{s_1,off}\phantom{:}\ &=0<t_\mathrm{s_2,off}<t_\mathrm{s_1,on}<t_\mathrm{s_2,on}<T\mbox{,\quad where}\\
t_\mathrm{s_1,off}:&\ \mbox{$y_2$ crosses $2\alpha_2/G-y_{0,1}$ from below to above,}&&\mbox{($u_1:1\to0$),}\\
t_\mathrm{s_2,off}:&\ \mbox{$y_1$ crosses $y_{0,2}$ from above to below,}&&\mbox{($u_2:1\to0$),}\\
t_\mathrm{s_1,on}:&\ \mbox{$y_2$ crosses $2\alpha_2/G-y_{0,1}$ from above to below,}&&\mbox{($u_1:0\to1$),}\\
t_\mathrm{s_2,on}:&\ \mbox{$y_1$ crosses $y_{0,2}$ from below to above,}&&\mbox{($u_2:0\to1$),}\\
T:&\quad\mbox{period of oscillation.}
\end{aligned}
\end{equation}

\paragraph{Alpha-type oscillations in the singular limit} For $G>G_1$ the periodic orbits consist of two coupled pairs of orbit segments, one pair of segments for pyramidal activity $y_1$  governed by the piecewise linear ODE \eqref{y1dot} (switching at times $t_\mathrm{s_1,off}$ and $t_\mathrm{s_1,on})$, and one pair of segments for inhibition $y_2$ governed by the piecewise linear ODE \eqref{y2dot} (switching at times $t_\mathrm{s_2,off}$ and $t_\mathrm{s_2,on}$), as shown in Figure~\ref{figsamllepstime}a. The coupling occurs through the switching times, as the switching times $t_\mathrm{s_2,on/off}$ are determined by $y_1$ crossing the threshold $y_{0,2}$, and the switching times $t_\mathrm{s_1,on/off}$ are determined by $y_2$ crossing the threshold $y_\mathrm{3,eq}-y_{0,1}=2\alpha_2/G-y_{0,1}$. Both piecewise linear ODEs \eqref{y1dot} and \eqref{y2dot} are affine ODEs of the form
\begin{align}
 \label{gen:ODE:affine}
  \ddot y&=b^2c-2b\dot y-b^2y, 
\end{align}
with inhomogeneity $c$.  The general solution $(y(t),\dot y(t))\in\mathbb{R}^2$ of \eqref{gen:ODE:affine} for its IVP with initial condition $(y(0),\dot y(0))=(y_0,\dot y_0)\in\mathbb{R}^2$ is
\begin{align}
    Y_\mathrm{adv}\left(t;b,c,\begin{bmatrix}
        y_0\\\dot y_0
    \end{bmatrix}\right)=\begin{bmatrix}
        y(t)\\ \dot y(t)
    \end{bmatrix}:&=\e^{M_bt}
    \begin{bmatrix}
        y_0\\\dot y_0
    \end{bmatrix}+\left[I_2-\e^{M_bt}\right]\e_1c\mbox{,\quad where\ }\nonumber\\
        \e^{M_bt}&=\e^{-bt}\begin{bmatrix}
       1+bt&t\\ -b^2t&1-bt
    \end{bmatrix}\mbox{,}\nonumber
\end{align}
$\e^{(\cdot)}$ is the exponential and $\e_1=(1,0)^\tran$ is the first standard basis vector. From this general formula we find that requiring periodicity of the pyramidal activity $y_1$ and the inhibition $y_2$ determines their initial conditions at the ``off'' switching time. The initial conditions
\begin{align*}
    \begin{bmatrix}
        y_1(t_\mathrm{s_1,off})\\
        \dot y_1(t_\mathrm{s_1,off})
    \end{bmatrix}&=\frac{2}{G}
    Y_\mathrm{per}(T,t_\mathrm{s_1,on}-t_\mathrm{s_1,off};1),\\
    \begin{bmatrix}
        y_2(t_\mathrm{s_2,off})\\
        \dot y_2(t_\mathrm{s_2,off})
    \end{bmatrix}&=\frac{2\alpha_4}{b^*}
    Y_\mathrm{per}(T,t_\mathrm{s_2,on}-t_\mathrm{s_2,off};b^*)\mbox{,\ where $Y_\mathrm{per}$ is defined as}\\
        Y_\mathrm{per}(T,\delta_t;b)&=\left[I_2-\e^{M_bT}\right]^{-1}\left[I_2-\e^{M_b(T-\delta_t)}\right]\e_1\mbox{}
\end{align*}
will lead to a solution of \eqref{y1dot}, \eqref{y2dot} that is periodic with period $T$. The switching times $t_\mathrm{s_1,on}$, $t_\mathrm{s_2,on/off}$ and the period $T$ satisfy $4$ relations that follow from the crossing conditions $y_1( t_\mathrm{s_2,off})=y_1( t_\mathrm{s_2,on})=y_{02}$ and $y_2( t_\mathrm{s_1,off})=y_2( t_\mathrm{s_1,on})=y_3-y_{01}$:
\begin{align}\label{pworbit:crossing}
\setlength{\arraycolsep}{1pt}
\begin{array}{rcllccl}
     y_{0,2}&=&\e_1^\tran Y_\mathrm{adv}(t_\mathrm{s_2,off}; &1,0, (\cdot))&\circ&\frac{2}{G}&Y_\mathrm{per}(T,t_\mathrm{s_1,on}\!-t_\mathrm{s_1,off};1)\\[0.5ex]
     y_{0,2}&=&\e_1^\tran Y_\mathrm{adv}(t_\mathrm{s_2,on}-T;&1,2/G, (\cdot))&\circ&\frac{2}{G}&Y_\mathrm{per}(T,t_\mathrm{s_1,on}\!-t_\mathrm{s_1,off};1)\\[0.5ex]
     \frac{2\alpha_2}{G}-y_{0,1}&=&\e_1^\tran Y_\mathrm{adv}(t_\mathrm{s_1,on}-t_\mathrm{s_2,off};&b^*,0,(\cdot))&\circ&\frac{2\alpha_4}{b^*}&
    Y_\mathrm{per}(T,t_\mathrm{s_2,on}\!-t_\mathrm{s_2,off};b^*)\\[0.5ex]
\frac{2\alpha_2}{G}-y_{0,1}&=&\e_1^\tran Y_\mathrm{adv}(-t_\mathrm{s_2,off};&b^*,2\alpha_4/b^*,(\cdot))&\circ&\frac{2\alpha_4}{b^*}&
    Y_\mathrm{per}(T,t_\mathrm{s_2,on}\!-t_\mathrm{s_2,off};b^*).
\end{array}
\end{align}
In \eqref{pworbit:crossing} and below we use the composition symbol $\circ$ with the meaning $f\circ x=f(x)$, $f\circ g\circ x=f(g(x))$ etc.\ to reduce nesting of brackets. Recall that we set $t_\mathrm{s_1,off}=0$ such that this switching time has been omitted in the defining equations \eqref{pworbit:crossing}. We have already inserted the appropriate initial conditions to ensure periodicity with period $T$ into \eqref{pworbit:crossing}, such that the system \eqref{pworbit:crossing} is a system of $4$ algebraic equations, depending on the $4$ unknown times $t_\mathrm{s_1,on}$, $t_\mathrm{s_2,off}$, $t_\mathrm{s_2,on}$ and $T$, and system parameters such  as $G$,  $b^*$ and the switching thresholds $y_{0,1},y_{0,2}$. For given switching times and period, one may obtain the complete periodic orbit as
\begin{align*}
    \begin{bmatrix}
        y_1(t)\\\dot y_1(t)
    \end{bmatrix}&=
        Y_\mathrm{adv}\left(t; 1,\frac{2}{G}u_1(t),(\cdot)\right)\circ \frac{2}{G}Y_\mathrm{per}(T,t_\mathrm{s_1,on}-t_\mathrm{s_1,off};1)\\
        \intertext{for $t\in[t_\mathrm{s_1,on}-T,t_\mathrm{s_1,on}]$}
    \begin{bmatrix}
        y_2(t)\\\dot y_2(t)
    \end{bmatrix}&=
        Y_\mathrm{adv}\left(t-t_\mathrm{s_2,off}; b^*,\frac{2\alpha_4}{b^*}u_2(t),(\cdot)\right)\circ  \frac{2\alpha_4}{b^*}Y_\mathrm{per}(T,t_\mathrm{s_2,on}-t_\mathrm{s_2,off};b^*)\\
        \intertext{for $t\in[t_\mathrm{s_2,on}-T,t_\mathrm{s_2,on}]$}
        \begin{bmatrix}
        y_3(t)\\\dot y_3(t)
    \end{bmatrix}&\!=\!    \begin{bmatrix}
        2\alpha_2/G\\0
    \end{bmatrix}\mbox{\quad where}\\\
         u_1(t)&=\begin{cases}
         0&\mbox{if $t\in[0,t_\mathrm{s_1,on}]$}\\
        1&\mbox{if $t\in[t_\mathrm{s_1,on}-T,0]$},
    \end{cases}\ 
    u_2(t)=\begin{cases}
         0&\mbox{if $t\in[t_\mathrm{s_2,off},t_\mathrm{s_2,on}]$}\\
        1&\mbox{if $t\in[t_\mathrm{s_2,on}-T,t_\mathrm{s_2,off}]$.}
    \end{cases}
\end{align*}
Figure~\ref{figsamllepstime}a shows an example of such a piecewise exponential and its switching times. The associated thresholds $y_\mathrm{3,eq}-y_{0,1}=2\alpha_2/G-y_{0,1}$ and $y_{0,2}$ for $y_2(t)$ and $y_1(t)$  for the switching events of $y_1(t)$ and $y_2(t)$ are shown as magenta and yellow horizontal lines.

\begin{figure}[ht] 
    \centering
    \includegraphics[width=0.9\textwidth]{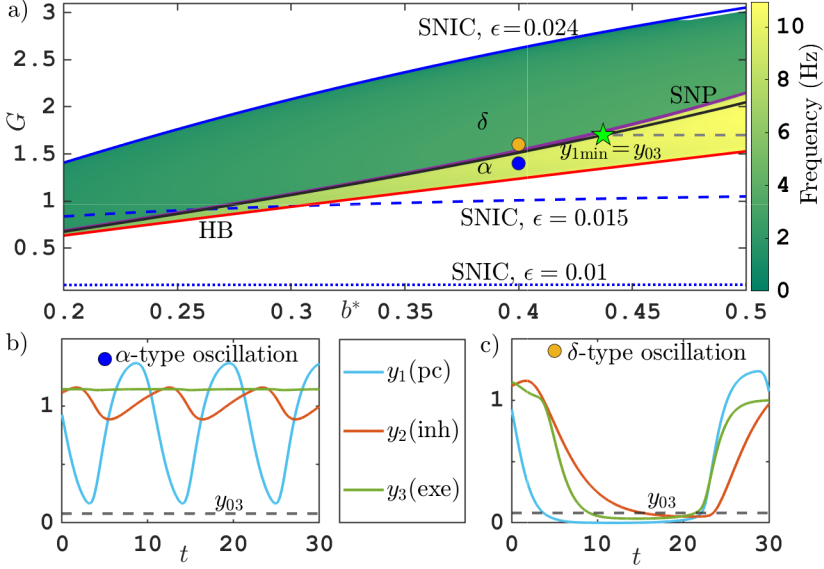}
    \caption{Panel (a): Composite bifurcation diagram overlaying regions of alpha and delta activity of  system~\eqref{nondimwitheps} for $\epsilon=0.024$, and the grazing bifurcation (black curve) detected by solving algebraic equations \eqref{pworbit:crossing},\,\eqref{pworbit:y1mindef} in the $(b^*,G)$ parameter plane.  
    Panel (b,c): time profiles of typical alpha activity ($(b^*,G)=(0.4,1.4)$, blue circle in panel a) and delta activity ($(b^*,G)=(0.4,1.6)$, yellow circle in panel a) for each neuron population at $\epsilon=0.024$. Folds of periodic orbit (SNP, two purple curves on top of each other) and Hopf bifurcation (HB, red) are for $\epsilon=0.024$, SNIC bifurcations are shown for $\epsilon=0.024,0.015,0.01$. Horizontal grey dashed line for $b^*\in[0.44,0.5]$ at $G=1.7$ is parameter path for continuation to grazing bifurcation reached at point labelled by green star. See Table \ref{tab:nondimparameters} for other parameters.}
    \label{fig5}
\end{figure}

\paragraph{Detection of grazing boundary for alpha activity}
The orbits given by algebraic system \eqref{pworbit:crossing} are alpha-type oscillations with constant support by excitatory interneuron population $y_3$. A criterion for the existence of these oscillations is that $y_1$ never crosses below the threshold $y_{0,3}$ at which the excitatory interneurons $y_3$ switch off. The piecewise exponential has a distinct minimum $y_{1,\min}$ at some time $t_{1,\min}$, labelled by a star in Figure~\ref{figsamllepstime}a. This minimum consistently occurs shortly after $t_\mathrm{s_1,on}$, in the interval between $t_\mathrm{s_1,on}$ and $t_\mathrm{s_2,on}$, due to the inertia caused by the second-order nature of the ODE for $y_1$. Hence, we may introduce $y_{1,\min}$ and $t_{1,\min}$ as additional parameters satisfying the relations
\begin{align}
    \begin{bmatrix}
        y_{1,\min}\\0
    \end{bmatrix}&=\begin{multlined}[t]
    Y_\mathrm{adv}(t_{1,\min}-t_\mathrm{s_1,on};1,2/G,(\cdot))\circ Y_\mathrm{adv}(t_\mathrm{s_1,on};1,0,(\cdot))\\
    \circ (2/G)Y_\mathrm{per}(T,t_\mathrm{s_1,on}-t_\mathrm{s_1,off};1).
    \end{multlined}
    \label{pworbit:y1mindef}
\end{align}
We add the two parameters $t_{1,\min}$ and $y_{1,\min}$ while tracking solutions of \eqref{pworbit:crossing},\,\eqref{pworbit:y1mindef} in system parameters to monitor the difference $y_{1,\min}-y_{0,3}$ for roots. The difference equals the distance between the point on the graph in Figure~\ref{figsamllepstime}a labelled with a star and the grey horizontal line $D_3=\{y_1=y_{0,3}\}$. When this difference is zero we have a grazing bifurcation: the alpha activity oscillation's orbit touches the threshold $D_3=\{y: y_1=y_{0,3}\}$ in a quadratic tangency. To find an orbit with such a quadratic tangency we fix $G=1.7$ and vary the time scale ratio $b^*$ from $b^*=0.5$ downwards along the horizontal dashed line shown in Figure~\ref{fig5}a. While doing so we monitor the quantity $y_{1,\min}-y_{0,3}$ until it reaches zero, which is the point that represents the grazing bifurcation. The grazing occurs at $(b^*,G)\approx(0.44,1.7)$, labelled with a green star in Figure~\ref{fig5}a.
\paragraph{Mapping alpha and delta oscillations: a two-parameter bifurcation analysis}
Following detection of this grazing bifurcation at $(b^*,G)\approx(0.44,1.7)$ we track its locus in two system parameters, time scale ratio $b^*$ and the feedback ratio $G$, by imposing the algebraic condition $y_{1,\min}-y_{0,3}=0$. The resulting curve in the $(b^*,G)$-plane is shown in Figure~\ref{fig5}a in black. We combine this curve with the other relevant bifurcations bounding alpha and delta activity:
\begin{description}
    \item[\textbf{(Hopf bifurcation)}] The Hopf bifurcation curve is the lower bound for alpha activity in $G$. It is given by \eqref{def:g1g3}  in the limit $\epsilon=0$ ($G_1=\alpha_2 b^*/(\alpha_4+y_{0,1}b^*/2)$), and shown in red in Figure~\ref{fig5}a for $\epsilon=0.024$. The value $\epsilon=0.024$ is the original value of $\epsilon$ given in Table~\ref{tab:nondimparameters}. Visually the Hopf curve for $\epsilon=0.024$ is indistinguishable from its $\epsilon=0$ limit.
    \item[\textbf{(Grazing bifurcation)}] Asymptotically for $\epsilon\to0$ the transition boundary between delta and alpha activity is the grazing curve determined by the algebraic equations \eqref{pworbit:crossing}, \eqref{pworbit:y1mindef}, shown in black in Figure~\ref{fig5}a.
    \item[\textbf{(SNIC bifurcation)}] The upper bound in $G$ for $\delta$ activity is a SNIC bifurcation, determined by the saddle-node of equilibria for the stable small-activity equilibrium, labelled $y_{\mathrm{eq}_1}$, and the small-activity saddle $y_{\mathrm{eq}_2}$ in Table~\ref{tab:xeqresults}. Table~\ref{tab:xeqresults} indicates that both equilibria exist for all positive values of $G$ for $\epsilon=0$ and (uniformly in $\epsilon$) positive thresholds $y_{0,1}$ and $y_{0,3}$. Hence, the parameter value $G_\mathrm{SNIC}$ for the SNIC bifurcation shown in Figure~\ref{dimandnondim}c has the limit $G_\mathrm{SNIC}\to0$ for $\epsilon\to0$. We included the SNIC bifurcation curve for three values of $\epsilon$ (the original value $0.024$ given in Table~\ref{tab:nondimparameters}, $0.015$ and $0.01$) in Figure~\ref{fig5}a to visualise the effect of decreasing $\epsilon$. For $(b^*,G)$ above the SNIC bifurcation curve the activity drops to its small-activity equilibrium.
\end{description}
To demonstrate the quantitative accuracy of the prediction  from the singular limit $\epsilon=0$ for the transition between delta and alpha activity, we include shading according to frequency in areas where we find oscillations for the positive value $\epsilon=0.024$ from Table~\ref{tab:nondimparameters}. We observe a sharp drop in frequency for increasing $G$ near the singular-limit grazing bifurcation (from yellow to green in Figure~\ref{fig5}a). Panels b and c in Figure~\ref{fig5} show time profiles at representative points in the alpha region ($(b^*,G)=(0.4,0.14)$, blue  circle) and in the delta region ($(b^*,G)=(0.4,0.16)$, yellow circle). For larger $b^*$ a small deviation is noticeable such that the effect of non-zero $\epsilon$ is slightly larger near $b^*\approx0.5$. The pair of folds of periodic orbits which coincides with the frequency drop (see Figures \ref{dimandnondim}b and \ref{dimandnondim}c) is also included for $\epsilon=0.024$ as a pair of purple curves, labelled SNP (the two SNP curves are close to each other, such that one obscures the other). The pair of folds overlaps with the grazing bifurcation for $\epsilon=0$ over large parts of the parameter plane.
% \hmq{in papare  they show this type of transition without fold of periodic orbit}

\section{Discussion and further analysis}
\label{sec:discussion}
We have identified the mechanism responsible for the transition between alpha- and delta-type activity in the Jansen-Rit model \citep{jansen1995electroencephalogram}. The transition between alpha-type oscillations and delta-type oscillations in neural circuits is relevant to applications because  alpha frequencies are associated with the promotion of Long-Term Potentiation (LTP), enhancing synaptic strength and facilitating memory formation during active learning, while delta frequencies are associated with Long-Term Depression (LTD), which promotes synaptic weakening and neural plasticity during deep sleep \citep{ks2021long}. The combined effects of these processes contribute to of learning via consolidation of memory across the wakefulness and sleep cycles \citep{burgdorf2023prefrontal,Brudzynski2024}. Therefore, understanding these transitions enhances our knowledge of how neural circuits regulate brain states and may provide insights into disruptions that affect sleep and cognitive functions.
In the Jansen-Rit model the transition is not associated with one of the generic codimension-one bifurcations, but it occurs over a small range of parameter values for decreasing excitatory feedback strength $A$ and with little to no hysteresis. In the transition the time profile and the frequency change sharply: oscillations during alpha activity show a shape typical for a single negative feedback loop because of the near-constant support from the excitatory interneurons. Oscillations during delta activity are of relaxation type with long periods of near-zero activity (see Figures~\ref{dimandnondim}d,e in Section \ref{sec:observation}). The feedback strength parameter $A$ affects how neurons respond to incoming signals by scaling the excitable population's output of firing rates \citep{jansen1995electroencephalogram}.

The results presented in Section \ref{sec:singpert:JR} show that the singular limit where one assumes that the activation is an all-or-nothing switch reduces the sharp transition in shape and frequency to a grazing bifurcation: when the pyramidal-cell activity ``touches'' the threshold line for switching off the excitatory interneurons in a quadratic tangency, the excitatory interneuron activity collapses. This causes further drop of pyramidal-cell activity, such that both populations approach near-zero activity. They stay there until the (slower) inhibition activity has dropped sufficiently low to permit recovery of excitatory activity in pyramidal and interneuron populations. 
%The result  obtained by numerical  bifurcation analysis and introducing singular limit to the  Jansen-Rit model by introducing small parameter that transforms the output firing rate into Heavside function 
(i.e., neuron response to incoming signals acts as binary operation off/on) in Section \ref{sec:singpert:JR}. 
The singular limit allow us to detect the curve of grazing bifurcation separating the two distinct wave rhythms (i.e., the boundary between two distinct types of periodic orbits)  in the two parameter plane $(b^*,G)$, where $b^*$ is the decay rate ratio of inhibitory to excitatory populations and $G$ is the ratio of postsynaptic amplitudes of inhibitory to excitatory populations in non-dimensionalised model \eqref{nondimwitheps}.

There are several further open questions for the analysis of the singular limit.

\paragraph{Small-threshold analysis of alpha-delta transition and SNIC}
Our analysis shows that the limit $\epsilon\to0$ with fixed thresholds $y_{0,1}$, $y_{0,3}$ does not support delta-type  oscillations. In delta-type  oscillations all potentials come close to zero such that these oscillations are only possible if there are no stable equilibria at small $y_1$, $y_3$. However, Figure~\ref{fig5} shows that for all sufficiently small $\epsilon$ there are no parameter values $(b^*,G)$ for which grazing has happened (so, above the black curve) but no equilibria with $y_1\approx 0$, $y_3\approx0$ exist (below blue solid/dashed/dotted curve called SNIC). This is because for positive $y_{0,1}$, $y_{0,3}$ and small $y_1$, $y_2$, $y_3$ the existing stable ``off'' equilibrium $y_{\mathrm{eq}_1}$ from Table~\ref{tab:xeqresults} attracts trajectories with $y_1,y_2,y_3\approx0$. If we also consider the thresholds for the excitatory populations as small quantities then delta-type oscillations are possible for small $\epsilon$: as we showed in Section~\ref{subsec:Equilibria in the  limit}, if the thresholds are proportional to $\epsilon$ then there are no small equilibria present for small $\epsilon$. A more detailed analysis is needed to find for which asymptotic order of magnitude for the thresholds a saddle-node (and, hence, SNIC) is present for parameters $G$ of order $1$.

The expressions in \eqref{pworbit:crossing} show that the dependence of the alpha-type periodic orbits on the threshold $y_{0,1}$ is regular, such that a replacing $y_{0,1}$ by zero in \eqref{pworbit:crossing} will not have a large effect. However, the grazing bifurcation is determined by the minimum of $y_1$ reaching $y_{0,3}$, such that for $y_{0,3}\to0$ the grazing bifurcation will shift. The pyramidal-cell potential $y_1$ decays exponentially as long as the inhibition $y_2$ is above the equilibrium value of the threshold $y_3-y_{0,1}\approx 2\alpha_2/G$. 
%\hmq{I think Jan here try to explain the effects of $b*$ we expect decreasing the decay rate of inhibition $y_2$ then this effects  the period of periodic exponential decay  for $y_1$ as we decrease $\epsilon$}
As a smaller $b^*$ leads to a slow-down in the change of inhibition $y_2$, a decrease of $b^*$ leads to a longer period of exponential decay for $y_1$, so to a smaller minimum of $y_1$. Hence, we expect the value of $b^*$ for the grazing bifurcation to decrease slowly with $\epsilon$, on the order of $1/|\log\epsilon|$, because of the exponential decay of $y_1$ for above threshold inhibition.

\paragraph{Canard explosion of alpha-type orbits at the singular Hopf bifurcation} The Hopf bifurcation, shown in the zoomed-in bifurcation diagram Figure~\ref{figbdsamlleps} can be studied in the singular limit $\epsilon=0$.
%\jsq{We may need the bifurcation diagram also for $Y_1$ or $y_1$ versus $G$, similar to $y_3$ vs $G$ in Fig.~\ref{dimandnondim}d for this. I suggest including the diagram for $\epsilon=0.024$ next to Figure~\ref{figbdsamlleps}.}
\begin{figure}[ht] 
    \centering   \includegraphics[width=\textwidth]{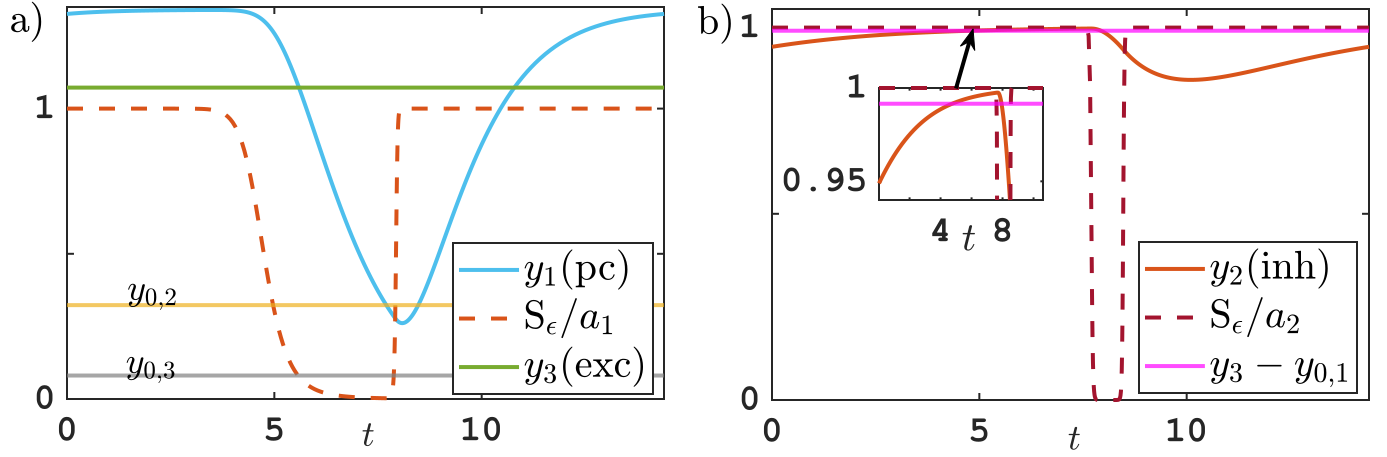}
    \caption{Time profile for alpha-type canard orbit near Hopf bifurcation (labeled HB1, in Fig~\ref{figbdsamlleps}) and activations 
$\newS_{\epsilon/a1}$ for $y_1$, and  $\newS_{\epsilon/a2}$ for $y_2$. Parameters: $G=1.49\approx G_1=1.48$ (green dashed vertical line in Fig~\ref{figbdsamlleps} near HB1), $\epsilon=0.001,y_{0,1},y_{0,3}=0.08$
and $y_{02}=0.3$, for others see Table \ref{tab:nondimparameters}. Inset in panel b shows that $y_2$ crosses the threshold $y_\mathrm{3,eq}-y_{0,1}$ of $y_1$.} 
    \label{canard}
\end{figure}
Figure~\ref{canard} shows that the periodic orbits in this canard-like explosion are parametrised by the time during which inhibition is switched off (drop of dashed curve in Figure~\ref{canard}b to zero), which is short for values of the parameter $G$ close to the bifurcation value $G_1$. Since the periodic orbits are piecewise exponentials, the asymptotic parameter-amplitude dependence  for $G\approx G_1$ can be determined explicitly.

\paragraph{Improved limit for slow inhibition} The small thresholds $y_{0,1}$ and $y_{0,3}$ and sharp sigmoid slopes in the activation function $\newS_\epsilon$ model the fact that excitatory populations have faster internal dynamics. We kept the threshold for the inhibition $y_{0,2}$ non-small, but replaced the activation function of the inhibition $\newS_{\epsilon/a_2}$ also by an all-or-nothing switch in the limit $\epsilon\to0$. A more appropriate limit consideration would be to keep a finite-slope sigmoid activation function $\newS_2(y_1-y_{0,2})$ in place for the inhibition and take the small-$\epsilon$ limit only for the excitatory populations. We would still encounter the same dynamical phenomenon for the alpha-delta transition in the limit, namely a grazing bifurcation. However, the limiting periodic orbits of alpha type would no longer be determined by an algebraic system with exponentials, such as \eqref{pworbit:crossing}, but as solutions of a nonlinear ODE over a finite time interval. So, finding the grazing bifurcation would require solving a nonlinear ODE numerically, instead of using the formulas in \eqref{pworbit:crossing}, \eqref{pworbit:y1mindef}.

\paragraph{Dependence of delta-type oscillations on parameters including non-zero input} Delta-type oscillations in the zero-input case studied here have excitation population potentials $y_1$ and $y_3$ at zero until the inhibition $y_2$ drops sufficiently low such that $y_1$ and $y_3$ can recover and generate a large-amplitude burst (see Figure~\ref{dimandnondim}e). The precise time ratios between time spend at zero and bursts depends on time scale ratios such as $b^*$, and thresholds. A positive external input $p$ will affect the time profile and frequency of delta-type oscillations, as positive $p$ will prevent $y_3$ from staying near zero, resulting in a positive slope at low potential values, as seen in the simulations of \citet{forrester2020role} and \citet{ahmadizadeh2018bifurcation}.

\section{Availability of Data} Scripts reproducing the computational data for all figures can be accessed at \url{https://github.com/jansieber/MSTA-alphadelta-anziam24-resources}
\section{Acknowledgments}
KTA gratefully acknowledges the financial support of the Engineering and Physical Sciences Research Council (EPSRC) via grant EP/T017856/1.
For the purpose of open access, the corresponding author has applied a ‘Creative Commons Attribution' (CC BY) licence to any Author Accepted Manuscript version arising from this submission.

Funding statement: Huda Mahdi gratefully
acknowledges the Higher Committee for Education Development (HCED),
Government of Iraq, for funding her PhD studies.
%\newpage
%\bibliographystyle{acm}
%\bibliography{references}
\bibliographystyle{abbrvnat}
%bibliographystyle{Vancouver}
\bibliography{reffirstyear}
\end{document}